\documentclass[10pt,oneside]{article}

\usepackage[utf8]{inputenc} %
\usepackage[T1]{fontenc}    %
\usepackage{url}            %
\usepackage{booktabs}       %
\usepackage{amsfonts}       %
\usepackage{nicefrac}       %
\usepackage{microtype}      %

\usepackage{lmodern}
\usepackage{times}

\usepackage{amssymb,amsmath,amsthm}
\usepackage{bbold}
\usepackage[short]{optidef}

\usepackage{framed,multirow,multicol}

\usepackage{xcolor,xspace}
\usepackage{lscape}
\usepackage{subcaption,graphicx,epsfig,tikz,caption}
\usepackage[normalem]{ulem}
\usepackage{enumerate}
\usepackage{verbatim}
\usepackage{makeidx,latexsym}
\usepackage[colorlinks=true,citecolor=blue]{hyperref}%

\usepackage{bm}
\usepackage{stmaryrd}
\usepackage{lscape,rotating}

\usepackage[capitalize]{cleveref}
\crefrangeformat{equation}{(#3#1#4)--(#5#2#6)}

\usepackage[numbers,sort&compress,square,comma]{natbib}

\usepackage{algorithmic,algorithm}

\usepackage[english]{babel}
\usepackage{bbm}
\usepackage{pgfplots}

    \usepackage{parskip}
\usepackage[letterpaper,margin=1.1in]{geometry}

    \usepackage{thmtools,thm-restate}

    \declaretheoremstyle[qed=$\square$]{definitionwithend}
    \declaretheorem[style=definitionwithend]{definition}
    
    \declaretheorem[style=definitionwithend]{example}
    
\PassOptionsToPackage{numbers, compress}{natbib}

\DeclareMathOperator{\dist}{dist}

\DeclareMathOperator{\conv}{conv}

\newcommand{\vx}{\mathbf{x}} %

\newcommand{\bbE}{\mathbb{E}}

\newcommand{\bbP}{\mathbb{P}}

\newcommand{\bbR}{\mathbb{R}}

\newcommand{\bbZ}{\mathbb{Z}}

\newcommand{\E}{\mathbb{E}}

\newcommand{\R}{\mathbb{R}}

\newcommand{\cF}{\mathcal{F}}

\newcommand{\cM}{\mathcal{M}}

\newcommand{\cP}{\mathcal{P}}

\newcommand{\cS}{\mathcal{S}}

\newcommand{\cX}{\mathcal{X}}
\newcommand{\cY}{\mathcal{Y}}

\newcommand{\st}{\text{s.t.}}

\newcommand {\beqn}{\begin{equation}}\newcommand {\eeqn}{\end{equation}}
\newcommand {\beqan}{\begin{eqnarray}}\newcommand {\eeqan}{\end{eqnarray}}
\newcommand {\beqa}{\begin{eqnarray*}}\newcommand {\eeqa}{\end{eqnarray*}}
\newcommand{\skipit}  [1] {}
\newcommand{\ignore}  [1] {}

\newtheorem{theorem}{Theorem}%

\newcommand{\cvar}{\operatorname{CVaR}}
\newcommand{\var}{\operatorname{VaR}}

\providecommand{\revise}[1]{{\protect{#1}}}

\begin{document}
\title{ Chance-Constrained Optimization \revise{under Limited Distributional Information: A Review of Reformulations Based on Sampling and Distributional Robustness}}
\author{Simge K\"{u}\c{c}\"{u}kyavuz\thanks{Industrial Engineering and Management Sciences, Northwestern University, Evanston, IL 60208, USA, \url{simge@northwestern.edu}} \and
        Ruiwei Jiang   \thanks{Industrial and Operations Engineering,  University of Michigan, Ann Arbor, MI 48109, USA, \url{ruiwei@umich.edu}}     
}

\date{\today}

\maketitle

\begin{abstract}
Chance-constrained programming (CCP) is one of the most difficult classes of optimization problems that has attracted the attention of researchers since the 1950s. 
In this survey, \revise{we focus on cases when only a limited information on the distribution is available, such as a sample from the distribution, or the moments of the distribution.}  We first review recent developments in mixed-integer linear  formulations of chance-constrained programs that arise from finite discrete distributions (or sample average approximation). We highlight successful reformulations and decomposition techniques that enable the solution of large-scale instances. We then review  active research in distributionally robust CCP, which is a framework to address the  ambiguity in the distribution of the random data. The focal point of our review is \revise{on} scalable formulations that can be readily implemented with state-of-the-art optimization software.  
Furthermore, we highlight the prevalence of CCPs with a review of applications across multiple domains. 
\end{abstract}

\section{Introduction} \label{sec:intro}

Most optimization models in practice involve problem parameters that are uncertain. Furthermore, in some cases these uncertain parameters involve risky outcomes with low probability. Therefore, requiring  feasibility of a solution  for every possible outcome may lead to overly conservative solutions. To remedy this, chance-constrained programming (CCP)  has emerged as a powerful paradigm to model system failure/reliability considerations and to  address the conservatism of a solution given a certain tolerance for risky outcomes.  

CCP is used to model risk-averse decision-making problems in a plethora of applications.
\revise{A few recent and active application domains---not meant to be an exhaustive list---include {\bf finance}, {\bf healthcare}, {\bf power systems}, {\bf transportation and routing}, {\bf supply chain, logistics, scheduling}, and {\bf wireless communications}.}
For example, in power systems, production levels need to be determined so as to meet peak load (demand) \cite[]{H93}. This problem is complicated by uncertainties in both generator availability (especially with renewables) and loads. The utility company's aim is to minimize the expected cost of power production while ensuring that the loss-of-load probability (i.e., the probability that the available generator capacity is insufficient to meet the peak load) is below an acceptable reliability level \cite[]{NKY90}.
In supply chain problems,  service level constraints are introduced  to limit the probability of stock-outs   \cite[]{chen2001}. In portfolio optimization problems, there is interest to restrict the downside risk at a certain threshold (value-at-risk) \cite[]{DJdVY08}.  Finally, in communications network design problems, a certain quality of service (QoS) with respect to packet losses needs to be ensured \cite[]{MR00}. Such risk, service, or reliability constraints are modeled using CCPs.  \revise{We review various applications of CCPs in Appendix~\ref{sec:app}. }

\subsection{Problem Definition} \label{sec:def}

Formally, for a given probability space $(\Omega, \cF, \bbP^0)$,  a chance-constrained  program (CCP) is given by 
\begin{subequations} \label{eq:ccp-orig}
\begin{align}
\min_{x} \quad & c^\top x\notag\\ 
\text{s.t.}\quad & \bbP^0(x\in \mathcal P(\omega) ) \geq 1-\epsilon,\label{eq:ccp-orig-cons}\\
& x \in \mathcal{X},
\end{align}
\end{subequations}
where $c\in \bbR^n$ is a cost vector, $\cX \subset \bbR^n$ represents a compact set defined by deterministic constraints on the decision variables $x$, possibly including integrality restrictions on some variables, $\omega \in \Omega \subset \bbR^d$ is a random vector with a  distribution $\bbP^0$, for a given $\omega$, $\mathcal P(\omega)$ represents the set of solutions that are safe or desirable,  and $\epsilon \in (0,1)$ is the risk tolerance for the decision vector $x$ being unsafe. For risk-averse
decision makers typical choices for the risk level are small
values, e.g.,  $\epsilon \le 0.05$. 
In this survey, we mainly focus on \emph{linear} chance constraints, i.e., polyhedral $\mathcal P(\omega)$. More precisely, let  
\begin{equation}\label{eq:polyhedral}
\mathcal P(\omega) :=\{x: T(\omega)x\ge r(\omega)\},
\end{equation}
where $T(\omega)$ is an $m\times n$ matrix of random constraint coefficients, and $r(\omega)\in \bbR^{m}$ is a vector of random right-hand sides.

Next, we introduce the taxonomy of CCPs. Constraint \eqref{eq:ccp-orig-cons} is said to be an \emph{individual} chance constraint for $m=1$, and a \emph{joint} chance constraint for $m>1$. If, for all $\omega\in \Omega$, we have  $T(\omega)=T$ for some deterministic $m\times n$ matrix $T$, and only  $r(\omega)$ is random, we say that the CCP has \emph{right-hand side (RHS) uncertainty.} In contrast, if the so-called \emph{technology matrix} $T(\omega)$ is random, we say that the CCP has \emph{left-hand side (LHS) uncertainty}, regardless of whether $r(\omega)$ is a fixed vector or is random. 
Most of the  work in CCP  can be seen as \emph{single-stage} (i.e., static) decision-making problems where the decisions are made here and now, and there are no  recourse actions once the uncertainty is revealed. In Section~\ref{sec:two-stage-ccp}, we discuss extensions to \emph{two-stage CCPs}. Finally, in many problems of interest, the decision vector $x$ is pure binary and this structure can be exploited to obtain stronger formulations and specialized algorithms. We refer to  such CCPs with pure binary variables as \emph{chance-constrained combinatorial optimization} problems.

CCP dates back to the early work of \citet{CCS58,CC63,MW65,prekopa-1970-on-proba}, and \citet{P73}, who first consider problems with individual   or joint chance constraints.  We refer the reader to  \cite{BL97,KW94,Dentcheva06,P95,prekopa2003ccp,shapiro-2009-lectur} for textbook treatment and detailed reviews that describe the earlier developments in this area. This survey is aimed at reviewing the developments in the past two decades primarily from a \revise{perspective of mathematical programming} reformulations for CCPs \revise{that can be easily implemented by practitioners using a state-of-the-art optimization solver, such as Gurobi and CPLEX. Beyond immediate use with reformulations, such solvers also facilitate the solution of  mixed-integer (conic) problems, using user-defined cutting plane or Benders decomposition approaches. We focus on the cases when there is limited knowledge of the probability distributions governing the random parameters under minimal assumptions (e.g., non-i.i.d.), and consider approaches based on sampling or distributional robustness. In Appendix \ref{sec:surveys}, we give Table \ref{tab:survey} that summarizes other reviews covering topics that are complementary.}

Despite long-standing interest and ubiquity in practice,  CCP remains one of the most challenging class of problems in general. 
There are two main challenges with CCPs.
\begin{enumerate}
  
    \item {\bf Non-convexity of the feasible set.} 
For certain special cases such as joint CCPs with RHS uncertainty  involving quasi-concave or log-concave distributions \cite{prekopa-1970-on-proba,P95,Wets1983,WETS1989}, or individual chance constraints with LHS uncertainty under a certain log-concave distribution and choice of $\epsilon$ \cite{Lagoa2005}, such as normal \cite{Kataoka1963}, there is an equivalent convex representation of the corresponding CCP. (See, also, \cite{vanAckooij-2011-survey} for a discussion on convexity.) 
In general, however, chance constraints even in the case with continuous $x$, polyhedral $\cP$, and only RHS uncertainty give rise to non-convex feasible regions in their original variable space.
We illustrate this challenge with an example. 

\begin{example}\label{ex:sen-ex}\cite{kuecuekyavuz-2012-mixin-sets}
Let $\omega_1$ and $\omega_2$ be dependent random variables with joint probability density function given in Table \ref{tab:pdf-ex}. 
Consider the CCP with RHS uncertainty
\begin{eqnarray*}
\min & x_1+x_2&\\
\mbox{s.t.}  & \bbP^0 \left \{
\begin{array}{lcl}
2x_1-x_2&\ge  &\omega_1\\
x_1+2x_2& \ge &\omega_2\\
\end{array} \right \}&\ge 0.6\\
&x\ge 0.&
\end{eqnarray*}
The   feasible region of this problem is non-convex as illustrated in  Figure \ref{fig:sen-ex}. \revise{It can be represented as a union of polyhedra whose extreme points are given by the so-called $(1-\epsilon)$-efficient points, which results in a feasible set that is not necessarily convex.}

\begin{table}[htb]
\caption{Joint probability density function of $\omega$}\centering \label{tab:pdf-ex}
\begin{tabular}{c|ccccccccc}
Scenario&1&2&3&4&5&6&7&8&9\\\hline
$\omega_1$ &0.75&0.5&0.5&0.25&0.25&0.25&0&0&0\\
$\omega_2$& 1.25&1.5&1.25&1.75&1.5&1.25&2&1.5&1.25\\
Probability  & 0.2 & 0.14& 0.06 & 0.06& 0.06& 0.3& 0.04 & 0.04 & 0.1
\end{tabular}
\end{table}

\begin{figure}[hbt]
\begin{center}
\resizebox{!}{4.5cm}{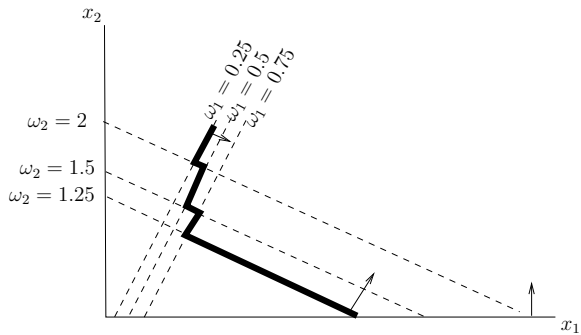}
\end{center}
\caption{The feasible region of the example CCP \revise{(from \cite{kuecuekyavuz-2012-mixin-sets})}.} \label{fig:sen-ex}
\end{figure}
\end{example}

Indeed, the resulting problems are NP-hard, in general \cite{luedtke-2008-integ-progr,nemirovski-2007-convex-approx}. \citet{BN98,calafiore-2004-uncer-convex-progr,CC06}, and \citet{NS05,nemirovski-2007-convex-approx}  approximate the non-convex chance constraint  with convex constraints such that the solution to this approximation is feasible with high probability. However, such methods could yield highly conservative solutions \cite[]{AS08} (see Section \ref{sec:approx}). 

  \item {\bf Difficulty of evaluating the probability of an undesirable solution.}  In \revise{general}, the distribution $\bbP^0$ in the chance constraint \revise{may not be} fully specified. A black-box simulation model or an oracle may be available to evaluate $\bbP^0$ for a given solution $x$, however it is not straightforward to integrate such  an oracle within the optimization model and the number of feasible solutions to evaluate is typically huge \cite{wu-2019-probab-partial}.
    In this survey, we focus on two main approaches to address this difficulty, namely the Sample Average Approximation (SAA) approach \cite[]{AS08} (Section \ref{sec:saa}) and the distributionally robust approach (Section \ref{sec:dro}).

\end{enumerate}

 There has been renewed and growing interest in CCP since the early 2000s \cite{DPR00,R02} to tackle these challenges. Capitalizing on the enormous success of mixed-integer programming (MIP) and conic optimization solvers since the early 2000s, our focal point is on reformulations that aim to circumvent the aforementioned challenges and enable progress towards the solution of this difficult class of problems.

\subsection{Outline}
Our survey is organized as follows. In the first part of this survey, in Section~\ref{sec:saa}, we consider CCPs under a finite discrete distribution. We consider a natural MIP formulation and valid inequalities for both RHS and LHS uncertainty in Sections \ref{sec:saa-rhs} and \ref{sec:saa-joint}, respectively. In Section \ref{sec:saa-alt-form}, we review alternative formulations and specialized methods for CCPs under a finite distribution.  In Section \ref{sec:two-stage-ccp}, we  describe a two-stage CCP and a Benders decomposition method for its solution. In Section \ref{sec:approx} we describe approximations of CCPs. 
In the second part of this survey, in Section \ref{sec:dro}, we consider distributionally robust CCPs, primarily under two types of uncertainty sets: moment-based (Section \ref{sec:moment}) and Wasserstein ambiguity sets (Section \ref{sec:wass}) and conclude in Section \ref{sec:conclude}.  We give some preliminaries, provide more details on some valid inequalities, and an overview of a wide range of applications in the Appendix.

\section{CCPs under Finite Discrete Distributions} \label{sec:saa}
In this section, we consider CCPs under a finite discrete probability space $(\Omega, 2^{\Omega}, \bbP_N)$, where $\Omega=\{\omega_1,\dots,\omega_N\}$ and $p_i=\bbP_N(\omega=\omega_i)$. Let $F_\omega$ denote the cumulative distribution function of $\omega$.  Of particular interest are such CCPs that result from 
the SAA approach  \cite{luedtke-2008-sampl-approx,Pagnoncelli2009}, which   approximates $\bbP^0$ via  a finite empirical distribution, $\bbP_N$. \revise{We refer the reader to Appendix \ref{sec:prelim} for preliminaries on two types of risk definitions for univariate random variables that are closely related to chance constraints and their formulations when the random variable follows a finite discrete distribution.}

  For ease of exposition, we will assume that the samples are 
 independent and identically distributed (i.i.d.) and consider the SAA formulation of CCP (i.e., $p_i=\frac{1}{N}, i\in[N]$). Throughout, for $a\in\bbZ_+$, let $[a]:=\{1,\dots,a\}$. The methods we discuss can be adapted to the case of non-i.i.d.\ scenarios, for example those that are obtained via importance sampling \cite{Barrera2016}. 
 
 The SAA formulation of \eqref{eq:ccp-orig} is
 \begin{subequations}\label{eq:saa-ccp}
\begin{align}
\min\limits_{x} \quad & c^\top x\\\ 
\text{s.t.}\quad & \frac{1}{N} \sum_{i \in [N]} \bm{1}(x\not \in \mathcal P(\omega_i) ) \leq \epsilon,  \label{eq:saa-ccp-cons}\\
& x \in \mathcal{X},
\end{align}
\end{subequations}
where $\bm{1}(\cdot)$ is the indicator function. From this formulation, it is apparent that the use of  finite discrete distribution circumvents the  difficulty of evaluating \revise{the violation probability when the distribution is not explicitly given}. Under non-equal probability scenarios, constraint \eqref{eq:saa-ccp-cons} is simply   $$\sum_{i \in [N]}p_i \bm{1}(x\not \in \mathcal P(\omega_i) ) \leq \epsilon.$$
When  $\cP(\cdot)$ is polyhedral as given by \eqref{eq:polyhedral}, 
formulation \eqref{eq:saa-ccp} for CCP under a  discrete distribution lends itself to an equivalent mixed-integer \revise{(linear)} program  via the introduction of binary variables and \revise{so-called} big-M constraints, \revise{which we will specify in the next section}. Hence, the non-convex feasible region in the original space of variables can be represented as \revise{a polynomial-sized} MIP with additional binary variables. \revise{While the resulting feasible region is non-convex, the MIP formulation  enables} the immediate use of   off-the-shelf MIP solvers. Next we present such MIP formulations for the RHS and LHS uncertainty cases.

\subsection{RHS uncertainty} \label{sec:saa-rhs}

First, let us consider the problem with RHS uncertainty. 
In this setting, the joint linear CCP   \eqref{eq:saa-ccp} with RHS uncertainty is reformulated  as a mixed-integer linear program  \cite{R02}  
\begin{subequations}\label{eq:ccp-re}
\begin{align}
\min\limits_{x,t,z} \quad& c^\top x\label{ccp-re-obj}\\
\st\quad& x\in\cX,\qquad Tx=\bar r+t, \label{eq:ccp-Tx}\\
&t_j\geq r_{i,j}(1-z_i),\quad\forall i\in[N],~\forall j\in[m],\label{ccp-bigM}\\
&\frac{1}{N}\sum_{i\in[N]}z_i\leq \epsilon,\label{ccp-knapsack}\\
&t\in\R_+^m,~z\in\{0,1\}^N\label{ccp-binary},
\end{align}
\end{subequations}
where $\bar r\in\R^m$ is chosen vector satisfying $r(\omega_i)\geq \bar r$ for all $i$ and ${r_i}=(r_{i,1},\ldots,r_{i,m})^\top$ denotes $r(\omega_i)-\bar r$. The choice of 
$\bar r$ ensures that the  vector $r_i$ is nonnegative for all $i\in[N]$. \revise{Since $\epsilon<1$, we have 
 $Tx\geq r(\omega_i)$ for some $i\in [N]$. Therefore, $Tx \geq \bar r$, and from \eqref{eq:ccp-Tx}, we have $t\in\R_+^m$  as stated in \eqref{ccp-binary}.} 
The binary variable $z_i$ encodes the indicator function in \eqref{eq:saa-ccp-cons} to  model the event $Tx\geq r(\omega_i)$. In particular,  if $z_i=0$, then  constraints~\eqref{ccp-bigM} enforce that $t\geq r_i$ holds and thus $T x\geq r(\omega_i)$ is satisfied. Otherwise, $z_i=1$,  and constraints~\eqref{ccp-bigM} reduce to the trivial relation $t\geq 0$. \revise{Such constraints are referred to as \emph{big-M constraints}, where in \eqref{ccp-bigM}, the big-M coefficients are given by $r_{i,j}$.} Finally, \eqref{ccp-knapsack} enforces that the probability of $x\not \in \cP(\omega)$  is within the risk threshold~$\epsilon$. Note that this constraint is equivalent to a cardinality constraint on the binary variables $\sum_{i\in[N]}z_i\leq \lfloor \epsilon N\rfloor =:k$. In the non-equiprobable case, it is a knapsack constraint $\sum_{i\in[N]}p_iz_i\leq \epsilon.$

In  the case of individual chance constraints, when $m=1$, we can linearize the 
\revise{chance constraint \eqref{eq:ccp-orig-cons}, where $\mathcal P(\omega)$ in \eqref{eq:polyhedral} is defined by a single constraint,} as $T x\ge F_\omega^{-1}(1-\epsilon)$ to lower bound the LHS with the $(1-\epsilon)$-quantile (see Definition \ref{def:var} in Appendix \ref{sec:prelim}). Under RHS uncertainty,  problems with joint chance constraints $(m>1)$ are more challenging. 
In fact, \citet{luedtke-2008-integ-progr} show that \revise{problem \eqref{eq:ccp-re}} is NP-hard  for $m>1$.  
\ignore{Constraints~\eqref{ccp-bigM}
are referred to as  big-$M$ constraints.} Often, formulations with big-M constraints result in weak LP relaxation bounds, which  hinder the convergence of the branch-and-bound  methods. 
Therefore, MIP approaches have focused on obtaining strong formulations for the SAA formulation to scale up the problem sizes that can be solved. To this end, an important substructure in the formulation~\eqref{eq:ccp-re} is given by the constraints \eqref{ccp-bigM} and \eqref{ccp-binary} for a fixed $j$. This particular substructure is a special case of the {\it mixing set} studied in  \cite{guenluek-2001-mixin-mixed} that involves general integer variables \revise{whose convex hull of solutions is given by the so-called mixing inequalities}. \revise{The specific form of the mixing inequalities}  involving only binary variables---known as the star inequalities---is independently considered in \citet{atamturk2000mixed} in the context of  vertex covering. \revise{In the ensuing discussion, we will specify the star (or binary mixing) inequalities with an additional strengthening in the context of CCPs.}

 We first consider strengthening based on an individual inequality in the \revise{joint} chance constraint. 
 More precisely, consider \eqref{ccp-bigM} and \eqref{ccp-binary} for a fixed $j$. We will drop the dependence on $j$ for notational convenience. 
The resulting system is nothing but a mixing set with binary variables given by
\[
\cM:=\left\{(t,z)\in\R_+\times\{0,1\}^N:~ t+r_{i}z_i \geq r_i,~\forall i\in[N]\right\}. 
\]
The (binary) mixing set $\cM$ involves $N$ inequalities   that share a  common continuous  variable $t$, but independent binary variables $z_i$, $i\in[N]$. 
The so-called  {\it mixing inequalities} of \citet{guenluek-2001-mixin-mixed} specialized to binary case, which is known to be equivalent to the so-called star inequalities introduced in~\cite{atamturk2000mixed},  are an exponential family of linear inequalities that provide the complete linear description of $\conv(\cM)$  (see also,  \citet[Theorem 18]{PW94}). Furthermore, this class of inequalities can be separated in polynomial time \cite{guenluek-2001-mixin-mixed,atamtuerk-2003-facet-mixed}, hence formulation \eqref{eq:ccp-re} can be strengthened using the mixing inequalities within a branch-and-cut framework.  
Somewhat surprisingly,  \citet{Kilinc-Karzan2019joint-submod} uncover that \revise{the  mixing inequalities for the set  $\cM$ are nothing but the   polymatroid inequalities valid for the epigraph of a submodular function $g(1-z)$ (in this case $g(1-z)=\max_{i\in [N]}\{r_i(1-z_i)\}$)} as defined in \citet{lovasz1983submodular}, \citet[Proposition 1]{atamturk2008polymatroids}. \revise{(See Appendix \ref{sec:submod} for a review of submodularity and polymatroid inequalities.)}

 \citet{luedtke-2008-integ-progr} further strengthen  formulation~\eqref{eq:ccp-re}  by exploiting 
the cardinality  constraint~\eqref{ccp-knapsack} and by studying 
 the resulting set given by \eqref{ccp-bigM}--\eqref{ccp-binary} for a fixed $j$. 
 In this case, an immediate strengthening is that of the big-M. Consider the \revise{mixing set with a cardinality constraint} 
 \[
\cM_C:=\left\{(t,z)\in\R_+\times\{0,1\}^N:~ t+r_{i}z_i \geq r_i,~\forall i\in[N], \sum_{i\in[N]}z_i\leq k\right\}. 
\]
 Sort  the values $r_i$ for
$i\in [N]$,
to obtain a permutation $\sigma$ such that:
\[r_{\sigma_1} \geq{} r_{\sigma_2} \geq{} \cdots{} \geq{} r_{\sigma_N}.\]  
Now observe that due to the cardinality constraint $\sum_{i\in[N]}z_i\leq k$, we must have $t\ge r_{\sigma_{k+1}}$. Therefore, we deduce that 
\[
\cM_C=\left\{(t,z)\in\R_+\times\{0,1\}^N:~ t+(r_{i}- r_{\sigma_{k+1}}) z_i \geq r_i,~\forall i\in[N], \sum_{i\in[N]}z_i\leq k\right\}. 
\]
Note, here, that this is an immediate big-M coefficient strengthening\revise{---the coefficient of $z_i$ is reduced from $r_i$ to $(r_{i}- r_{\sigma_{k+1}})$---}that can be readily incorporated into the MIP \revise{formulation \eqref{eq:ccp-re}}. This strengthening uses the quantile information that $t\ge r_{\sigma_{k+1}}$.

Due to their common usage, we give a precise definition of the resulting mixing inequalities that make use of the cardinality-based strengthening next. 
 Consider a subset
$S = \left\lbrace s_1, s_2, \ldots{}, s_\ell
\right\rbrace \subseteq{} \left\lbrace\sigma_1, \sigma_2, \ldots{},
  \sigma_k\right\rbrace$ such that $ r_{s_i}\geq{} r_{s_{i+1}}$ for $i=1,
\ldots{}, \ell$, where $s_1 = \sigma_1$ and $s_{\ell+1} =
\sigma_{k+1}$.  \citet{luedtke-2008-integ-progr} show that a strong mixing inequality valid for $\cM_C$ is given by 
\begin{align}
t+ \sum_{i = 1}^{\ell} \left( r_{s_i} -  r_{s_{i+1}} \right)
z_{s_i} \geq{} r_{s_1},
\label{eq:mixing-gen}
\end{align}
\revise{which is precisely in the form of a basic mixing inequality proposed for a mixing set without the cardinality constraint, with the exception that $S\subseteq{} \left\lbrace\sigma_1, \sigma_2, \ldots{},
  \sigma_k\right\rbrace$ with $r_{s_{\ell+1}} =
r_{\sigma_{k+1}}$ as opposed to $S\subseteq N$ and $r_{s_{\ell+1}} =0$ for a basic mixing inequality.} 
\revise{\citet{kuecuekyavuz-2012-mixin-sets} show that} this idea can be adapted to the non-equiprobable case by  \revise{setting $\varphi:=\arg\max\{j:\sum_{i=1}^j p_{\pi_i} \le \epsilon\}$, where $\pi$ is a permutation of the scenarios in increasing order of their probabilities. This choice of $\varphi$ gives a valid cardinality restriction on the binary variables}.  Furthermore, inequality \eqref{eq:mixing-gen} can be  strengthened  by further use of the cardinality relation or for the case where the scenarios are not equiprobable when constraint \eqref{ccp-knapsack} is in the form of a knapsack inequality  \cite{luedtke-2008-integ-progr,kuecuekyavuz-2012-mixin-sets,abdi2016mixing-knapsack,zhao2017joint-knapsack}. \revise{We describe this further strengthening in inequality \eqref{eq:mixing-individual} in  Appendix \ref{sec:blending}.}

Next, we illustrate the strengthened mixing inequalities on our numerical example (Example \ref{ex:sen-ex}).  Consider the first inequality inside the chance constraint and note that $\varphi=3$ with respect to $\omega_1$. Note that the scenarios are already ordered in nonincreasing order with respect to the possible values of  $r_1(\omega)$.  Therefore, we have $t_1\ge 0.25=r_1(\omega_4)$. A possible strengthened mixing inequality is for $S=\{1,3\}$ given by 
$$t_1+(0.75-0.5)z_1+(0.5-0.25)z_3\ge 0.75.$$ It is easy to see the validity of this inequality. If $z_1=0$, then we must have $t_1\ge 0.75$, which satisfies this inequality. If $z_1=1$ and $z_3=0$, then we must have $t_1\ge 0.5$, which is also satisfied. Finally, when $z_1=z_3=1$, the inequality reduces to $t_1\ge 0.25$, which holds due to the $(1-\epsilon)$-quantile relation.

So far, we reviewed inequalities based on an individual inequality inside the \revise{joint} chance constraint. If we consider multiple inequalities inside the chance constraint jointly, the resulting set is an intersection of multiple mixing sets that share a common set of binary variables $z$, but independent continuous variables $t_j, j\in [m]$. For this case,  \citet[Theorem~3]{atamturk2000mixed} show that adding the mixing inequalities written for each set to the LP relaxation of the set defined by ~\eqref{ccp-bigM} and~\eqref{ccp-binary} is sufficient to  obtain the convex hull of solutions \revise{in the absence of cardinality/knapsack constraints}. Furthermore,  \citet{Kilinc-Karzan2019joint-submod} show how to extend their framework exploiting submodularity to recover this result, as well as extend it to propose the so-called \emph{aggregated mixing inequalities} that incorporate lower bounds on the continuous variables based on the quantile relation. 
For the special case of two-sided chance constraints, the valid inequalities in the convex hull description provided in \citet{liu2018intersection} are equivalent to the aggregated mixing inequalities. The aggregated mixing inequalities do not directly use the cardinality information, but use it indirectly through the lower bound on the continuous variables obtained from the quantile. In contrast,~\citet{kuecuekyavuz-2012-mixin-sets} and \citet{zhao2017joint-knapsack} propose valid inequalities for \revise{the intersection of mixing sets} by directly considering the cardinality/knapsack constraint. \revise{We elaborate on these approaches in Appendix~\ref{sec:blending}.}

\subsection{LHS uncertainty}\label{sec:saa-joint}

Now consider the problem with uncertainty data in both LHS and RHS. 
In this setting, the joint linear CCP   \eqref{eq:saa-ccp} with LHS uncertainty is reformulated  as a mixed-integer linear program  \cite{R02}  
\begin{subequations}\label{eq:ccp-re-lhs}
\begin{align}
\min\limits_{x,z} \quad& c^\top x\label{ccp-re-obj-lhs}\\
\st\quad& x\in\cX,\\
&T(\omega_i)x\geq r(\omega_i)-M(\omega_i)z_i,\quad\forall i\in[N],\label{ccp-bigM-lhs}\\
&\frac{1}{N}\sum_{i\in[N]}z_i\leq \epsilon,\label{ccp-knapsack-lhs}\\
&z\in\{0,1\}^N\label{ccp-binary-lhs},
\end{align}
\end{subequations}
where $M(\omega_i), i\in [N]$ is a vector of big-M coefficients such that when $z_i=1$, inequality \eqref{ccp-bigM-lhs} is redundant. 

In Section~\ref{sec:saa-rhs} we exploited the mixing structure associated with~\eqref{ccp-bigM} and~\eqref{ccp-binary} for a fixed $j$. In other words, we considered an individual inequality inside the (joint) chance constraint. Furthermore, we considered RHS uncertainty only.   In this section we will consider LHS as well as RHS uncertainty, and we will jointly consider the inequalities inside the chance constraints for any $m\ge 1$. 

The mixing \revise{inequalities} described in Section \ref{sec:saa-rhs} relies on the fact that all scenarios share the same LHS for a given $j\in [m]$, that is $t=T_j x \revise{-\bar r_j}$, where $T_j$ is the $j$th row of $T$. Due to this relation, we arrive at a mixing structure with $N$ constraints that share the same continuous variable $t$ and different binary variables. In contrast, in LHS uncertainty case, we no longer have a common continuous variable. Can we still apply the mixing \revise{inequalities}?

As it turns out, we can indeed extend the mixing \revise{inequalities} to generate other classes of  valid inequalities for joint chance-constrained programs with LHS uncertainty. To do so, we solve the following single-scenario optimization
problem
for all scenarios $\omega \in \Omega$ and for a given $\phi \in \bbR^n$:
\begin{subequations}\label{eq:g-mixing-subproblem}
\begin{align}
q_\omega\left( \phi \right) = \min\limits_{x}\quad & \phi^{\top} x \\
 &x\in \cP(\omega), & \\
& x\in \mathcal X. &
\end{align}
\end{subequations}
We  sort  the values $q_\omega\left( \phi{} \right)$ for
$\omega\in\Omega$,
to obtain a permutation $\sigma$ such that:
\[q_{\sigma_1}\left( \phi{} \right) \geq{} q_{\sigma_2}\left(
  \phi{} \right) \geq{} \cdots{} \geq{} q_{\sigma_{N}}\left( \phi
\right).\] Observe that $\phi^{\top} x \ge q_{\sigma_{k+1}}(\phi)$ is a valid inequality. Furthermore, substituting $t=\phi^{\top} x$ and $r=q(\phi)$ in inequality \eqref{eq:mixing-gen}, we obtain a valid inequality of the desired form.  These inequalities are referred to as \emph{quantile cuts}. This and related inequalities based on  quantile information have been studied in~\cite{ahmed2017nonanticipative,luedtke-2014-branc-and,song-2014-chanc-const,qiu2014covering-lp,LKL16,xie2018quantile}. These inequalities consider the interaction between the decision variables across multiple inequalities in the chance constraint, which results in improved  computational performance. In another line of work, \citet{Tanner2010} propose a class of cuts based on the irreducibly infeasible subsystems (IIS) of an LP that requires that a subset of scenarios are satisfied. The authors demonstrate the efficacy of this approach in a vaccine allocation application.

\ignore{consider a subset
$S = \left\lbrace s_1, s_2, \ldots{}, s_\ell
\right\rbrace \subseteq{} \left\lbrace\sigma_1, \sigma_2, \ldots{},
  \sigma_k\right\rbrace$ such that $ q_{s_i}\left( \phi{}
\right) \geq{} q_{s_{i+1}}\left( \phi{} \right)$ for $i=1,
\ldots{}, \ell$, where $s_1 = \sigma_1$ and $s_{\ell+1} =
\sigma_{k+1}$. Then it follows that the inequality
\begin{align}
\phi x+ \sum_{i = 1}^{\ell} \left( q_{s_i}\left( \phi{}
\right) -  q_{s_{i+1}} \left( \phi{} \right) \right)
z_{s_i} \geq{} q_{s_1}\left( \phi{} \right)
\label{cut:feasiblity_ideal-gen-q}
\end{align}
is a valid inequality.
}

\subsection{Alternative formulations and methods}\label{sec:saa-alt-form}
 While we focus on natural big-M formulations that can be easily adopted by practitioners, it is important to note that there are alternative reformulations for this class of problems relying on the concept of $(1 - \epsilon)$-efficient points, which
 are an exponential number of points representing the multivariate value-at-risk associated with the chance constraint \eqref{join-prob} to be specified later. \revise{We provide a brief overview here, and refer the reader to \cite{Dentcheva06}, and references therein, for a more detailed treatment of  methods based on $(1 - \epsilon)$-efficient points.}

\begin{definition}\cite[]{prekopa-1990-dual-method}\label{def:peff}
Let $\nu\in \bbR^m$  be such that $F_\omega(\nu)\ge 1-\epsilon$ and $F_\omega(\nu-\varepsilon)<1-\epsilon$ for $\varepsilon\ge \mathbf{0}$, $\varepsilon\ne \mathbf{0}$. The point $\nu$  is called {\it $(1-\epsilon)$-efficient}.
\end{definition}

In Example~\ref{ex:sen-ex}, observe that  $\nu\in\{(0.25,2), (0.5,1.5), (0.75,1.25)\}$  is $(1-\epsilon)$-efficient. The $(1-\epsilon)$-efficient points then prescribe the extreme points of the non-convex feasible region as seen in Figure~\ref{fig:sen-ex}.

There are several methods in the literature that rely on the enumeration of the exponentially many $(1-\epsilon)$-efficient points  \cite{prekopa-1990-dual-method,S92,DPR00,Kogan2014,Kogan2016}.   Such alternative formulations lead to specialized branch-and-bound algorithms described in
 \cite{beraldi-2002-probab-set,BR02,R02,saxena-2010-mip-refor}.  
\citet{S92} uses the $(1-\epsilon)$-efficient points to  give a disjunctive programming reformulation of joint chance constraints with finite discrete distributions. Valid inequalities are proposed based on the extreme points of the reverse polar of the disjunctive program, which can be separated by a cut generation linear program (CGLP) \cite[]{B79}. 
 \citet{kuecuekyavuz-2012-mixin-sets}  gives a compact and tight extended formulation based on disjunctive programming for $m=1$. \citet{Vielma2012} extend this formulation for varying $m>1$ to obtain a hierarchy of stronger relaxations.  
\citet{DPR00}  use $(1-\epsilon)$-efficient points to obtain various reformulations of probabilistic programs with discrete random variables, and to derive valid  bounds on the optimal objective function value.
\citet{R02} uses the concept of $(1-\epsilon)$-efficient points to derive consistent orders on different scenarios representing the discrete distribution. The consistent ordering is represented with precedence constraints, and valid inequalities for the resulting precedence-constrained knapsack set are proposed. \citet{BR02} propose a branch-and-bound method for probabilistic integer programs using a partial enumeration of the $(1-\epsilon)$-efficient points. \revise{In a related line of work, \citet{L12} introduces the concept of $(1-\epsilon)$-sufficient points---a concept less computationally demanding than $(1-\epsilon)$-efficient points as described in \cite{Lejeune2012AoR}---leading to a Boolean reformulation method and MIPs for CCPs. \citet{LejeuneMargot2016} extend this method to  quadratic CCPs. }

Alternatively, \citet{ahmed2017nonanticipative} consider a Lagrangian relaxation of the MIP formulation by creating copies of the variables, and relaxing the non-anticipativity constraint that these variables are equal. The authors derive extended formulations \revise{(without big-M coefficients)} whose relaxations achieve stronger bounds than the basic formulation (without mixing strengthening). \revise{In addition,~\citet{ahmed2017nonanticipative} propose a heuristic scheme to generate a conservative approximation for CCP and later~\citet{jiang2020alsox} show that this approximation is tighter than the classical CVaR approximation (see Definition~\ref{def-cvar} and  \citet{nemirovski-2007-convex-approx}).}

Furthermore, for problems with pure binary variables  and special  structures, i.e., for \emph{combinatorial CCPs}, stronger formulations 
have been  developed (see, e.g., \cite{beraldi-2009-exact-approac,song-2013-branc-and,song-2014-chanc-const,wu-2019-probab-partial,Hong2015,LK18}). 
For example,  \citet{song-2014-chanc-const} study chance-constrained bin packing problems, and propose a formulation that does not involve additional indicator variables to represent \eqref{eq:saa-ccp-cons} based on the so-called lifted probabilistic cover inequalities. Later, \citet{wang2019chance} consider a closely related formulation with multiple chance constraints and derive lifted cover, clique, and projection inequalities based on a bilinear reformulation. In a related line of work,~\citet{wang-2019-solut-approac} consider a chance-constrained assignment
problem and its distributionally robust variant, and propose  lifted cover
inequalities based on a bilinear reformulation of the problem. For chance-constrained knapsack problems, \citet{yoda2016convexity} provide sufficient conditions for the convexity of the formulation, \citet{klopfenstein2008robust}, \citet{de2018boolean}, \citet{han2016robust}, and~\citet{joung2020robust} derive approximate but more tractable formulations that can provide near-optimal solutions, and \citet{goyal2010ptas} derive a fully polynomial time approximation scheme when the random item sizes are independent and Gaussian. In addition, \citet{nikolova2010approximation} studies approximation algorithms for general chance-constrained combinatorial optimization problems with random parameters following either the Gaussian distribution or a general distribution. \citet{xie-2020-bicrit-approx} provide a bicriteria approximation algorithm for a class of chance-constrained covering problems and their distributionally robust variants that finds a solution within constant factor of the violation probability and a constant factor of  the optimal objective. 

For chance-constrained set covering models with RHS uncertainty,~\citet{beraldi-2002-probab-set,
  saxena-2010-mip-refor} propose a specialized branch-and-bound algorithm based on the enumeration of
$(1-\epsilon)$-efficient points. Subsequently,~\citet{saxena-2010-mip-refor} derive 
 polarity cuts to improve the
computational performance of this approach. For individual chance-constrained set-covering problems with LHS
uncertainty,~\cite{fischetti-2012-cuttin-plane} developed  cutting plane approaches 
for the case that all components of the Bernoulli random vector
\(\omega_i\) are independent. In
addition,~\citet{wu-2019-probab-partial} propose an exact approach for a partial set covering problem for the case that there exists an oracle to retrieve the probability of any events
under \(\mathbb{P}^0\). In another line of work, \citet{goyal2008approximation} and \citet{swamy2011risk} propose approximation algorithms for chance-constrained set-covering problems with optimality guarantees.

In addition to the aforementioned combinatorial CCPs, \citet{padberg1989branch} and \citet{campbell2008probabilistic} study chance-constrained traveling salesman problems, \citet{song2016risk} incorporate a chance constraint into a bi-level shortest path interdiction problem, and \citet{ishii1981stochastic} and \citet{geetha1993stochastic} study chance-constraint variants of the spanning tree problem.

It bears mentioning that there are recent nonlinear programming-based approaches that directly address the non-convexity of chance constraints. 
\citet{CAA06} give a global optimization algorithm that successively partitions the non-convex feasible region until a global optimal solution is obtained. \citet{TTN95} give an algebraic geometry algorithm for a scheduling problem with joint chance constraints that solves a series of chance-constrained integer programs with varying reliability levels.   \citet{Pena-Ordieres2020} derive smooth non-convex reformulations of the chance constrained based on the sampled empirical distribution. Other nonlinear programming approaches, which may result in solutions that are stationary points, include difference-of-convex optimization methods \cite{Hong2011}, sequential outer and inner approximations \cite{Geletu2017}, and sequential  cardinality-constrained quadratic optimization methods \cite{Curtis2018}.

Finally, throughout, we have assumed that the risk level $\epsilon$ is fixed. However, in practice, the decision-maker may be interested in the trade-offs between risk level and the optimal objective. One way to assess this would be to solve the problem for multiple values of fixed $\epsilon$. Alternatively, \citet{Shen2014}  proposes a novel \emph{variable risk threshold} model in which the risk tolerance is adjustable with an appropriate penalty function in the objective to prevent high risk.  The author proposes a MIP formulation for this problem for individual chance constraints. \citet[Theorem 8]{xie-2019-optim-bonfer} show that the corresponding optimization problem is strongly NP-hard.    \citet{ElciNB2018} propose a stronger MIP formulation for this problem under RHS uncertainty. Finally, \citet{Lejeune2016} consider joint chance constraints also with LHS uncertainty and propose a Boolean-based mathematical formulation for this model.

\subsection{Two-stage Chance-Constrained Programming}\label{sec:two-stage-ccp}

Thus far, we have considered a decision-making problem that is static. In other words, the decisions are made here-and-now before the revelation of the outcome of a random event. However, in most practical situations, there are multiple decision stages---intervened by a probabilistic event---and the decision-maker  takes recourse actions in the later epochs based on the observed outcome of the event. In this section, we focus on problems that involve two stages. For example, in a power generation setting, the day-ahead problem determines the on/off status of the conventional generators a day before realizing the demand (load) or supply (in case of renewable generators). Then the second-stage problem ensures that the loss-of-load probability is no more than a pre-specified risk level $\epsilon\in(0,1)$. Therefore, a two-stage chance-constrained model is called for.

As before, the random outcome  $\omega$ is defined on a probability space $(\Omega,2^{\Omega},\bbP_N)$. Let $\mathbb{E}[\cdot]$ denote the expectation operator taken with respect to $\omega$.
\citet{LKL16}  propose the  two-stage chance-constrained mixed-integer program
\begin{subequations}
\label{eq:cc-1stage}
\begin{align}
\min\limits_{x}  &\quad c^\top x + \bbP_N\left(x \in \mathcal P( \omega) \right) \mathbb E[h(x, \omega)|x\in \mathcal P( \omega)],\\
&\quad\bbP_N(x\in \mathcal P( \omega))\ge 1-\epsilon \label{join-prob}\\
& \quad x\in \mathcal X,
\end{align}
\end{subequations}
where $\mathcal P(\omega) =\{x:\exists y \text{ satisfying } W(\omega)y\ge r(\omega)-T(\omega)x, y\in \mathcal Y\}$  and the second-stage problem is given by 
\begin{subequations}\label{eq:cc-2stage}
\begin{align}
h(x,\omega) =\min\limits_{y} &  \quad g(\omega)^\top y\\
& \quad W(\omega)y\ge r(\omega)-T(\omega)x\\
& \quad y\in \mathcal Y.
\end{align}
\end{subequations}
Here, $g(\omega)$ is a vector of second-stage objective coefficients, $\cY$ is the domain of the second-stage decision vector $y$. For a  related model that considers only the feasibility of the second-stage problem  without an associated  second-stage cost function $h(x,\omega)$, we refer the reader to \cite{luedtke-2014-branc-and}.

The two-stage chance-constrained  problem can be formulated as a large-scale  mixed-integer program by introducing a big-$M$  term for each inequality in the chance constraint and a binary variable for each scenario. In particular, analogous to the static CCP, 
the deterministic equivalent formulation (DEF) of the two-stage CCP  may be stated as
\begin{subequations}\label{eq:cc-def}
 \begin{align}
&&\min\limits_{x,y,z}  & \quad c ^\top x+ \frac{1}{N}  \sum\limits_{i\in [N]}
g(\omega_i) ^\top y(\omega_i)z_i\label{eqn:DEF_obj}\\
&& &\quad T(\omega_i)x +W(\omega_i)y(\omega_i)\geq r(\omega_i)-M(\omega_i)z_i, & i \in [N] \label{eqn:DEF_const_3}\\
&& &\quad \frac{1}{N} \sum\limits_{i \in [N]}   z_i\leq \epsilon, \label{eq:gen-cc}\\
&& &\quad x\in \mathcal X, y(\omega_i)\in \mathcal Y,& i\in [N]\\
&&&\quad z_i\in\{0,1\} &i\in [N],
\end{align}
\end{subequations}
where $z_i, i\in [N]$ is a binary variable that equals 0 only if the second-stage problem for scenario $\omega_i$ has a feasible solution, and $M(\omega_i)$ is a vector of large enough constants that makes constraint \eqref{eqn:DEF_const_3} redundant if $z_i=1$, i.e., if the second-stage problem for scenario $\omega_i$ need not be feasible.
The rest of the constraints are interpreted similarly as before. 

This formulation  poses multiple challenges in addition to the usual difficulties of a formulation with big-M  constraints \eqref{eqn:DEF_const_3}. First, the objective function \eqref{eqn:DEF_obj} is nonlinear. Second,   the problem is large scale due to the copies of the variables $y(\omega_i)$  and  the   large number of binary variables $z_i$ for $i\in [N]$. Nevertheless, the formulation \eqref{eq:cc-def} has a decomposable structure---for a fixed first-stage vector $x$, the problem decomposes into independent scenario problems. Furthermore, if  $y$ is a continuous decision vector and $\cY$ is polyhedral, then the second-stage problems are linear programs. Next we describe a Benders-type decomposition  algorithm that not only exploits this decomposable structure, but also replaces the weak big-$M$  constraints \eqref{eqn:DEF_const_3}\ with stronger optimality and feasibility cuts, using the mixing structure.

\subsubsection{Benders Decomposition-Based Branch-and-Cut Algorithm}
\label{sec:solutionmethod}
Benders method  \cite{B62}, or its specific use in the classical two-stage stochastic programming  (without chance constraints) referred to as the $L$-shaped method \cite{VW69}, is the method of choice for \revise{two-stage stochastic programs, where the second-stage problems are linear programs. 
However, despite their similar structure to two-stage CCPs,} these methods are not immediately applicable to \eqref{eq:cc-def}, since  the feasibility and optimality cuts of the Benders method \revise{ensure that
 all second-stage problems are feasible and the second-stage costs of all scenarios are considered}, which is not the case for two-stage CCPs, \revise{where violations are allowed in the second-stage problems and the second-stage costs  of  the feasible second-stage problems are accounted for}. For general recourse problems, feasibility and optimality cuts different from the traditional Benders cuts must be developed. 

Let $\eta_i$ represent a
lower bounding approximation of the optimal objective function value of the
second-stage problem under scenario $\omega_i, i\in [N]$. Without loss of generality, we assume that  $\eta_i\ge 0, i\in [N]$. At each iteration of a Benders decomposition method, a sequence of
  relaxed master
problems (RMP) are solved:
\begin{subequations}\label{eq:cc-rmp}
\begin{align}
&& \min\limits_{x,z,\eta}\quad & c ^\top x+  \frac{1}{N}\sum\limits_{i\in [N]}\eta_i \label{eq:rmp-1st}\\
&& &\frac{1}{N}\sum\limits_{i \in [N]}   z_i\leq \epsilon, \\
&&&(x,z)\in \mathcal{F},&\\
& &&(x,z,\eta)\in \mathcal{O},&\\
&& &x\in \mathcal X \\
&&&z\in\{0,1\}^{N},\label{eq:rmp-last}
\end{align}
\end{subequations}
where, $\mathcal{F}$ and $\mathcal{O}$ denote the set defined by the feasibility
and optimality cuts---to be specified later---respectively.

At iteration $k$, let $(x^k,z^k)$ be the optimal solution to the RMP.  Given this first-stage solution, suppose that we solve the  LP \eqref{eq:cc-2stage} for outcome $\omega$ to obtain $h(x^k,\omega)$.
The feasibility cuts in set $\cF$ are derived from the solution to this LP. If $z^k_i=0$
for some $i \in [N]$, then the second-stage problem must be feasible.  If it is
infeasible for a scenario $j\in[N]$, then there exists an extreme ray $\psi_{ \omega_j}$ associated with the dual of
\eqref{eq:cc-2stage} for scenario $\omega_j$ that yields the inconsistent 
solution. Then, letting 
$\phi=\psi_{\omega_j}^\top T(\omega_j)$ in \eqref{eq:g-mixing-subproblem}   gives a violated valid \revise{mixing} inequality that cuts off this infeasible solution $(x^k,z^k)$.  
If, on the other hand, for all $\omega \in \Omega$, the second-stage problem associated with scenario $\omega$ such that $z^k(\omega)=0$ is indeed feasible, then the current solution $(x^k,z^k)$ is a feasible solution and no feasibility cuts are necessary. However, optimality cuts may be needed. Next we describe how to obtain valid optimality cuts.

Let $\psi_{\omega_j}$ be the dual vector associated with the
optimal basis of the  second-stage problem
\eqref{eq:cc-2stage} for scenario $\omega_j$ at this iteration. One possible big-M optimality cut is given by \cite{cloud,WGW12} 
  \begin{equation} \label{eq:bigM-optcut}
\eta_j + M_j z_j \geq  \psi_{\omega_j}^\top ( r(\omega_j) - T(\omega_j) x),
\end{equation}
where $M_j, j\in [N]$ is a big-M coefficient vector.

Next we describe a stronger optimality cut proposed by \cite{LKL16} that leads to faster convergence to an optimal solution.
Clearly, the traditional Benders optimality cut,  $\eta_j\geq \psi_{\omega_j} ^\top
(r(\omega_j)-T(\omega_j) x)$ is a valid optimality cut for
$x\in \mathcal X$ (in fact for $x\in \mathcal P(\omega) $) if $z_j=0$. However, it may not be
valid  for all $x\in\mathcal X$ for solutions with $z_j=1$. To obtain a valid optimality cut, 
we solve the following secondary problem with $\phi=\psi_{\omega_j}^\top
T(\omega_j) $:
\begin{align}
\bar v_{\omega_j}(\phi)=\min\limits_{x,y}  &\quad\phi x\notag\\  \quad & x \in  \mathcal X,\quad y\in \cY.\notag
\end{align}
Then we add the optimality cut of the
form \begin{equation}\eta_j  + \left(\psi_{\omega_j} ^\top r(\omega_j)-\bar
v_{\omega_j}(\phi)\right)z_j\geq \psi_{\omega_j} ^\top (r(\omega_j)-T(\omega_j) x)
.\label{ideal-cut-type1}\end{equation}
To see the validity of this inequality at $z_j=1$, note that in this case, the second-stage objective function contribution for scenario $\omega_j$ is zero. Furthermore,  inequality \eqref{ideal-cut-type1} evaluated at  $z_j=1$ reduces to $\eta_j\ge \bar
v_\omega(\phi) -\phi x $. Because $\bar
v_\omega(\phi) -\phi x \le 0$ for all $x\in  \mathcal X$ and $\eta_j\ge 0$, this inequality is trivially satisfied.
The finite convergence of the resulting algorithm is proven in \cite{LKL16} under certain assumptions.

In Table~\ref{tab:rhs}, we summarize a set of computational experiments that appear in \cite{LKL16} to show the effectiveness of the approaches discussed so far. The instances are based on a resource planning problem adapted from \cite{luedtke-2014-branc-and}.   In the first stage, the number of  servers among $s$ types of servers to employ is determined. The second-stage problem is to allocate the servers to clients of $\tau$ types, so that their demands are met with high probability ($1-\epsilon$). Instances with various choices of $N, \epsilon, \tau, s$ are tested and we report  the average statistics for three random instances  generated for the combination reported in each row. 
We compare the proposed ``Strong" decomposition algorithm which uses the optimality cuts  \eqref{ideal-cut-type1} with DEF \eqref{eq:cc-def} and the  decomposition approach (referred to as  ``Basic") which uses the mixing-based feasibility cuts and  the big-$M$ optimality cuts \eqref{eq:bigM-optcut} with an appropriate choice of big-$M$ as described in \cite{LKL16}. We  report the solution times (in seconds)  only for Strong decomposition, because for DEF and Basic, all instances tested reach the time limit of one hour. We also report the percentage optimality gap at termination under the Gap column. 
 In most cases, DEF is unable to find a feasible solution to the LP relaxation, as indicated by a `-'. In cases when it is able to find a feasible solution, it ends with a gap ranging from 4\% to 8\%. On the other hand, Basic is able to find a feasible solution for all instances, but is unable to prove optimality for any of the 36 instances tested.  It ends after an hour with optimality gaps ranging from 2\% to 7\%. In contrast, the Strong decomposition algorithm, based on the proposed strong optimality cuts, is able to solve most of the instances to optimality. For the two unsolved instances (indicated by a superscript $^1$ under the Gap column), the average optimality gap is  less than 0.1\%. These results highlight the importance of using strong formulations and decomposition for large-scale instances.

\begin{table}[ht]
\center
\caption{Result for instances with random RHS (summarized from \cite{LKL16})}
\label{tab:rhs}
\begin{tabular}{||c|c|c|c|c|c|c|c|c|}
\hline 
\multicolumn{2}{ ||c| }{Instances}     	  &\multicolumn{1}{ |c| }{DEF}    &\multicolumn{1}{ |c| }{Basic }    &\multicolumn{2}{ |c| }{Strong } \\    
 \hline				           
$(N, \epsilon)$   &$(s,\tau)$ 		& \texttt{Gap (\%) }    				 & \texttt{Gap  (\%) }  		 & \texttt{Time }  & \texttt{Gap  (\%) }   \\
\hline
\multirow{3}{*}{ (2000, 0.05)}     &(5,10)      &  4.60			& {2.34}		 				&166		&0	   			\\
     						   &(10,20)     &   -			 	& {2.93}     	 				&483		&0			\\	
     						   &(15,30)	&   -				& {2.69}             			&1106		&0				\\						   
 \hline
\multirow{3}{*}{ (2500, 0.05)}     &(5,10)      &   4.64			& 2.61    						&279		&0			\\
     						   &(10,20)     &   - 			& 3.08  						&711		&0	\\	
						   &(15,30)	&   -			   	& 2.88						&1819		&0.09$^1$   \\
 \hline

\multirow{3}{*}{ (2000, 0.1)}     &(5,10)      	&   7.1            	& 5.46    						&723		&0			\\
     						  &(10,20)    	&  - 				& 5.99     			 	&1069    	&0			\\	
						  &(15,30)    	&  -				& 6.27			   		   &1032		&0                 \\
 \hline  
 
\multirow{3}{*}{ (2500, 0.1)}     &(5,10)      	&  7.63      		& 5.32  					 	&641		&0			\\
     						  &(10,20)     	&   -  			& 5.79     				 	&1198 	&0		\\	
      					  &(15,30)     	&   -  			& 6.03    					 	&2112		&0.02$^1$		\\

 \hline

\end{tabular}
\end{table}

It is important to note that in this model, the undesirable outcomes $\omega$ such that $x\not\in \mathcal P(\omega) $ are simply ignored.  \citet{LKL16} propose an extension of the two-stage model \eqref{eq:cc-1stage}, where they allow so-called recovery decisions for the undesirable scenarios. They discuss how to resolve a potential time inconsistency in two-stage CCP. Furthermore, the Benders decomposition-based solution method is extended to operate in the case of recovery. 

\citet{Elci2018} extend this framework to a two-stage chance-constrained optimization model with a mean-risk objective, using the conditional value-at-risk as a risk measure. The authors apply this framework to a humanitarian relief network design problem and demonstrate its effectiveness on a case study based on hurricane preparedness in Southeastern United States. 
\citet{Lodi2019} extend this two-stage framework to  convex second-stage problems, motivated by hydro-power scheduling applications. They  build an outer approximation of the nonlinear second-stage formulations to design a Benders-type algorithm that converges to an optimal solution under mild assumptions. They demonstrate the computational benefit of the decomposition algorithm on a case study based on hydroplant data from Greece.

We close this subsection by noting the assumption of continuous second-stage variables can be lifted by leveraging the developments for decomposition algorithms for classical two-stage stochastic mixed-integer programs, where the second-stage problems also involve integer decisions \cite{CT1997,LL1993,NSen2008,Ntaimo2013,EORMS2010,NSen2005,SH2005,SS2006,GKS2013,ZK14,QS16,SMIPTutorial2016}. These methods rely on iteratively convexifying the second-stage problems and updating the feasibility and optimality cuts accordingly. These methods can be combined with the Benders-type algorithm we described to enable the solution of two-stage CCPs with integer variables at the second stage.

\subsection{Approximations}\label{sec:approx}

Given the difficulty of solving the exact formulations of CCPs or their SAA reformulations, one line of research has focused on inner and outer approximations of CCPs that are more tractable. This tractability often comes at the price of conservatism in the resulting solutions.  Here  we briefly review these formulations and refer the reader to \cite{Ahmed2018,Lejeune2018} for a review of relaxations and approximations for CCPs. 

\begin{itemize}
    \item {\bf Scenario approximation.} Scenario approximation (SA) \cite[e.g.,][]{calafiore-2004-uncer-convex-progr,CampiGaratti2008,CC06,campi,deFarias-vanRoy-2004} entails sampling to approximate the distribution $\bbP^0$ with a finite distribution $\bbP_N$ with a set of outcomes  $\Omega=\{\omega_1,\dots,\omega_N\}$. However, unlike the SAA model \eqref{eq:saa-ccp},  a usual stochastic program (not chance-constrained) is solved enforcing that the relations inside the chance constraint hold for each scenario. Thus, the scenario approximation problem is given by  
\begin{subequations} \label{eq:ccp-scen-app}
\begin{align}
\min_{x} \quad & c^\top x\notag\\ 
\text{s.t.}\quad & x\in \mathcal P(\omega), & \omega\in \Omega,\label{eq:ccp-scen-app-cons}\\
& x \in \mathcal{X},
\end{align}
\end{subequations}
As a result, for polyhedral $\cP(\omega)$ and continuous $x$, the resulting SA formulation is a large-scale LP. The authors give a finite sample  guarantee that the solution to this problem is feasible to the original CCP with high probability. Interestingly, this sample size does not depend on $m$, under certain assumptions.  Unfortunately, the required sample size is typically large and the resulting solution 
is overly conservative. The SAA approach \cite{luedtke-2008-sampl-approx,Pagnoncelli2009} is aimed at alleviating the conservatism of the SA approach by enforcing the chance constraint, with a smaller risk level, over the finite distribution $\bbP_N$, albeit as a  MIP  as opposed to an LP.

    \item {\bf CVaR approximation.} 
From the Definitions~\ref{def:var} and~\ref{def-cvar} of value-at-risk (VaR) and conditional value-at-risk (CVaR), respectively, it is readily apparent that for a univariate random variable $X$, $\cvar_{1-\epsilon}(X)\ge \var_{1-\epsilon}(X)$. Therefore, for individual chance constraints $(m=1)$, one can approximate the constraint $\bbP(r(\omega)-T(\omega)x\le 0)\ge 1-\epsilon$, or in other words, $\var_{1-\epsilon}(r(\omega)-T(\omega)x)\le 0$ with  $\cvar_{1-\epsilon}(r(\omega)-T(\omega)x))\le 0$. For the case of finite discrete distributions, this approximation leads to tractable reformulations due to the LP representation of CVaR given in \eqref{cvar_classic}. In particular, for individual chance constrained CCP \eqref{eq:saa-ccp}, the CVaR approximation LP is 
\begin{align*}
\min_{x} \quad & c^\top x\notag\\ 
\text{s.t.}\quad & \eta+\frac{1}{\epsilon N}\sum_{i\in [N]}w_i \le 0,\\
& w_i\geq r(\omega_i)-T(\omega_i)x - \eta,~\forall~i\in [N],\\
& x \in \mathcal{X}.
\end{align*}

In general, though, it is not possible to represent CVaR tractably  \cite{nemirovski-2007-convex-approx}. Nevertheless, \citet{nemirovski-2007-convex-approx} give a family of safe (i.e., feasible with high probability) and, in some cases, tractable approximations---referred to as generator-based approximations---that include the Bernstein approximation \cite{Pinter1989}. They show that the tightest such approximation is a CVaR approximation. However, CVaR approximation is also  conservative in some cases \cite{Alexander-var-cvar-2004}.  We refer the reader to \cite{NEMIROVSKI2012}, and references therein, for a survey on related safe tractable approximations for individual chance constraints.

In the case of joint chance constraints ($m>1$), it is worthwhile to note that even for the discrete case, while a vector-valued multivariate VaR definition exists (Definition~\ref{def:peff}), there is no unified definition of multivariate CVaR (see \cite{merakli-2018-vector-valued} and the discussions therein). This poses challenges in formulating related CVaR-based approximations that are tractable. One approach is to scalarize the multivariate random vector $r(\omega)-T(\omega)x$ and use the corresponding univariate CVaR.   Considering the ambiguity of the scalarization weights leads to a multivariate CVaR definition that can be represented as a challenging MIP with big-M constraints \cite{noyan-2013-optim-with}. MIP strengthening techniques can be used to improve the computational performance of the resulting multivariate CVaR formulations \cite{kuccukyavuz2016cut,liu-2017-robus-multic,Noyan2019}.

\item {\bf Bonferroni approximation.}  Given that joint chance constraints are significantly harder than individual chance constraints, one approximation scheme that is commonly considered replaces the joint chance constraint with $m$ individual chance constraints. In this case, consider replacing the joint chance constraint $\bbP(T_j(\omega)x\ge r_j(\omega), j\in [m])\ge 1-\epsilon$ with 
\begin{align}
&\bbP(T_j(\omega)x\ge r_j(\omega))\ge 1-\epsilon_j, \label{eq:bonf1} \\
&\mbox{where } \sum_{j\in [m]} \epsilon_j\le \epsilon. \label{eq:bonf2}
\end{align}
From Bonferroni's inequality, it follows that any solution satisfying 
constraints \eqref{eq:bonf1}--\eqref{eq:bonf2} also satisfies the joint chance constraint \cite{nemirovski-2007-convex-approx,ChenSimSun-robust-persp-2007}. Because optimizing over $\epsilon_j$ is, in general, difficult, a common choice is $\epsilon_j=\epsilon/m, j\in [m]$. However, this is also known to be a conservative approach \cite{chen-2010-from-cvar,nemirovski-2007-convex-approx}.

\end{itemize}

Note that while these approximations provide some statistical guarantees for feasibility, they are known to be conservative and do not come with optimality guarantees. Indeed, \citet{xie-2020-bicrit-approx} show an inapproximability result for CCPs.   \citet{Ahmed2014} uses a similar idea as \cite{nemirovski-2007-convex-approx}, this time to obtain a convex (Bernstein) relaxation that yield deterministic lower bounds.  Integrated chance constraints proposed by  \citet{KleinHaneveld1986}  replaces the non-convex chance constraints with a quantitative measure of shortfalls that lead to polyhedral representations \cite{Haneveld2006} in the discrete case. In this case, they are equivalent to the LP relaxation of the MIP formulation of CCP. Alternatively, statistical lower bounds can be obtained by using order statistics based on SAA solutions \cite{luedtke-2008-sampl-approx,Pagnoncelli2009}. 
 Such deterministic or statistical bounds are useful in assessing the  quality of a solution obtained from an approximation.

The finite sample guarantees of sampling based methods \cite{calafiore-2004-uncer-convex-progr,campi,CC06,luedtke-2008-sampl-approx,Pagnoncelli2009} 
are much too large and conservative in practice. On the other hand, for small $N$, the out-of-sample performance of the SAA solution may even be infeasible to the original problem. For example, in \cite{wu-2019-probab-partial}, the authors consider a partial set covering problem when an oracle that can evaluate the true probability of the desired event is available. They observe that for sample sizes that lend themselves to a tractable solution of the resulting MIP, the SAA solution is often infeasible to the original problem. This is related to the over-fitting phenomenon in machine learning when the solution of the problem is highly sensitive to the samples $\{\omega_i \}_{i \in [N]}$ used to obtain it. In the next section, we describe an approach that alleviates this problem.

\section{Distributionally Robust Chance-Constrained Programming} \label{sec:dro}

 Given the unavailability of the exact distribution $\bbP^0$ and the potential overfitting issues due to SAA-based approaches, there has been growing interest in modeling stochastic optimization problems that are \emph{distributionally robust} (see \cite{Rahimian2019DistributionallyRO} and references therein).  

Formally,  a \emph{distributionally robust chance-constrained program} (DRCCP) is modeled as
\begin{subequations}\label{eq:dr-ccp}
\begin{align}
\min_{x} \quad & c^\top x\\\ 
\text{s.t.}\quad & \sup_{\bbP \in \cF(\beta)} \bbP(x \not\in \cP(\omega)) \leq \epsilon \label{eq:dr-ccp-cons}\\
& x \in \mathcal{X},
\end{align}
\end{subequations}
where $\cF(\beta)$ is an ambiguity set of distributions and $\beta$ is a set of parameters that describe the ambiguity set. Accordingly, the distributionally robust chance constraint \eqref{eq:dr-ccp-cons} ensures that the chance constraint is satisfied with respect to all distributions in $\cF(\beta)$, even the worst possible one.

Several types of ambiguity sets have been studied in the literature based on various characteristics of the distribution, including moments, shape information (e.g., symmetry and unimodality), support, mixture models, and discrepancy measures (e.g., Wasserstein and $\phi$-divergence) ~\cite{ghaoui-2003-worst-case,erdogan-2005-ambig-chanc,calafiore-2006-distr-robus,hanasusanto-2015-distr-robus,vandenberghe-2007-gener-cheby,nemirovski-2007-convex-approx,zymler-2011-distr-robus,xu-2012-optim-under,ahmed-2013-probab-set,jiang-2015-data-driven,chen-2018-data-driven,yang-2018-wasser-distr,li-2019-ambig-risk,xie-2019-distr-robus,Lasserre2019}. These ambiguity sets lead to different computational tractability and conservatism of the corresponding DRCCP. In this survey, we will focus on moment-based ambiguity sets (Section \ref{sec:moment}) and Wasserstein ambiguity sets (Section \ref{sec:wass}).

\subsection{Moment-based ambiguity}\label{sec:moment}

There are many successful developments on the tractability of single and joint chance constraints with moment ambiguity sets, which characterize \(\mathcal{P}\) based on moment information of
\(\mathbb{P}^0\)~\cite{calafiore-2006-distr-robus,zymler-2011-distr-robus,yang-2014-distr-robus,hanasusanto-2015-distr-robus,xie-2016-deter-refor,hanasusanto-2017-ambig-joint,li-2019-ambig-risk}.

For known mean value $\mu$ and covariance matrix $\Sigma$,~\citet{ghaoui-2003-worst-case} characterize a moment ambiguity set
\begin{align*}
    \cF(\mu, \Sigma) := \{\bbP: \E[\omega] = \mu, \E[(\omega - \mu)(\omega - \mu)^{\top}] = \Sigma\}.
\end{align*}
All probability distributions in $\cF(\mu, \Sigma)$ need to have the designated first two moments, and are otherwise allowed to have different distribution types (e.g., Gaussian, Gaussian mixture, etc.) or different support (e.g., discrete or continuous). Perhaps surprisingly,~\citet[Theorem 1]{ghaoui-2003-worst-case} show that DRCCP is second-order conic representable for individual chance constraints (i.e., $m = 1$). Specifically, if $T(\omega) := \omega^{\top}A + T_0$ for some data matrix $A \in \bbR^{d \times n}$ and vector $T_0 \in \bbR^{1 \times n}$ and $r(\omega) := b^{\top}\omega + r_0$ for some data vector $b \in \bbR^{d}$ and constant $r_0 \in \bbR$, then constraint \eqref{eq:dr-ccp-cons} is equivalent to
\begin{align}
\mu^{\top} (b - Ax) + \sqrt{\frac{1 - \epsilon}{\epsilon}}\|\Sigma^{1/2}(b - Ax)\|_2 \leq T_0 x - r_0. \label{eq:icc-socp}
\end{align}
This indicates that DRCCP may improve not only the out-of-sample performance of CCP when the sample size $N$ is small but also the computational tractability. The same result is also discovered by~\citet{calafiore-2006-distr-robus} and~\citet{wagner2008stochastic}. In addition,~\citet{zymler-2011-distr-robus} point out an interesting fact that, for $m=1$ and ambiguity set $\cF(\mu, \Sigma)$, constraint \eqref{eq:dr-ccp-cons} is equivalent to its conservative approximation that replaces the chance constraint with CVaR (see Definition~\ref{def-cvar}), i.e., $\displaystyle \sup_{\bbP \in \cF(\mu, \Sigma)} \mbox{CVaR}_{1 - \epsilon}(r(\omega)-T(\omega)x  ) \leq 0$. 

For individual chance constraints, the result of~\citet{ghaoui-2003-worst-case} can be extended in multiple directions while maintaining both \emph{exactness} and \emph{computational tractability}. For example,~\citet{cheng2014distributionally} incorporate support information into $\cF(\mu, \Sigma)$ (e.g., specifying that $\bbP$ is supported on a convex set) and derive an exact reformulation of \eqref{eq:dr-ccp-cons} based on linear matrix inequalities. \citet{zhang-2018-ambig-chanc} consider potential errors of estimating the mean value $\mu$ and covariance matrix $\Sigma$, e.g., when this is done based on  inadequate historical data. To address this, they adopt an alternative ambiguity set proposed by~\citet{delage-2010-distr-robus} to allow the true mean value of $\omega$ to be within an ellipsoid centered at $\mu$ and the true covariance matrix to be bounded from above by $\Sigma$. For this extended ambiguity set,~\citet{zhang-2018-ambig-chanc} show that constraint \eqref{eq:dr-ccp-cons} is still second-order conic representable. For ambiguity set $\cF(\mu, \Sigma)$,~\citet{xu-2012-optim-under} study a distributionally robust variant of the stochastic dominance constraint (see, e.g.,~\citet{dentcheva2003optimization}), which requires different risk tolerances for violating a chance constraint with different magnitudes. More precisely, they study constraints $\displaystyle \sup_{\bbP \in \cF(\mu, \Sigma)} \bbP[T(\omega)x \geq r(\omega) - s] \leq \epsilon - {\beta}(s)$ for all $s \geq 0$, where ${\beta}(s)$ is a pre-specified non-decreasing function of $s$, and show that these constraints are conic representable for various ${\beta}(s)$ functions. Furthermore,~\citet{yang-2014-distr-robus} and~\citet{xie-2016-deter-refor} consider an extension that allows the event $x \in \cP(\omega)$ to depend non-linearly on $x$ and $\omega$, e.g., $x \in \cP(\omega)$ if and only if $f(x, \omega) \geq 0$, where function $f(x, \omega)$ is concave in $x$ and quasiconvex in $\omega$. For example,~\citet[Corollary 2]{yang-2014-distr-robus} recast \eqref{eq:dr-ccp-cons} as a linear matrix inequality if $r(\omega)$, as well as each entry of $T(\omega)$, is either convex quadratic or linear in $\omega$.

It is also possible to extend~\citet{ghaoui-2003-worst-case} by incorporating shape information into the ambiguity set $\cF(\mu, \Sigma)$. For example,~\citet[Lemma 3.1]{calafiore-2006-distr-robus} strengthens $\cF(\mu, \Sigma)$ by additionally requiring $\bbP$ to be \emph{centrally symmetric} (that is, $\bbP[A] = \bbP[-A]$ for any Borel set $A \subseteq \bbR^d$) and derives a conservative approximation of constraint \eqref{eq:dr-ccp-cons}.~\citet{hanasusanto2015decision} considers a similar ambiguity set and allows the true covariance matrix to be bounded from above by $\Sigma$ (instead of matching it exactly as in $\cF(\mu, \Sigma)$). Consequently,~\citet[Theorem 3.4.3]{hanasusanto2015decision} recasts \eqref{eq:dr-ccp-cons} as a set of conic constraints. Different from~\cite{calafiore-2006-distr-robus},~\citet[Theorem 1]{li-2019-ambig-risk} strengthens $\cF(\mu, \Sigma)$ by requiring that $\bbP$ is \emph{$\alpha$-unimodal} (a generalized notion of unimodality; see~\citet{dharmadhikari1988unimodality} for definition). They show that constraint \eqref{eq:dr-ccp-cons} is equivalent to a set of second-order conic constraints.~\citet[Example 3.4.4]{hanasusanto2015decision} considers a similar ambiguity set, which bound the true covariance matrix from above by $\Sigma$, and recasts \eqref{eq:dr-ccp-cons} as linear matrix inequalities.~\citet{stellato2014data} also considers a similar ambiguity set as in~\citet{li-2019-ambig-risk} but requires $\bbP$ to be centered around $\mu$. In that case,~\citet[Section 4.1.1]{stellato2014data} recasts \eqref{eq:dr-ccp-cons} as a single second-order conic constraint. There are works that consider other shape information and provide \emph{tractable conservative approximations} of \eqref{eq:dr-ccp-cons} (i.e., maintaining computational tractability at a potential cost of exactness). For example,~\citet{ChenSimSun-robust-persp-2007} replace the covariance information in $\cF(\mu, \sigma)$ with bounds on forward and backward deviations, which capture the \emph{asymmetry} of $\bbP$, and derive a conservative approximation of \eqref{eq:dr-ccp-cons} via second-order conic constraints.~\citet{li2018distributionally} drop the covariance restriction from $\cF(\mu, \Sigma)$ while adding in that $\bbP$ is \emph{log-concave} and supported on an ellipsoid centered at $\mu$. For this case,~\citet{li2018distributionally} derive conservative and relaxing approximations of \eqref{eq:dr-ccp-cons}, all via second-order conic constraints.~\citet{postek-2018-rbst-opt} replace the covariance information in $\cF(\mu, \Sigma)$ with the mean absolute deviation (MAD) from the mean and further require that $\omega$ is componentwise \emph{independent}. For that case,~\citet{postek-2018-rbst-opt} derive a conservative approximation of \eqref{eq:dr-ccp-cons} based on second-order conic constraints.

The special case of \emph{combinatorial DRCCPs} with individual chance constraints is in general \emph{intractable} because of the binary decision variables. Nevertheless, various formulation strengthening and algorithmic techniques can be applied to solve these problems more effectively. For example, \citet{ahmed-2013-probab-set} exploit supermodularity of their distributionally robust set covering problem to derive a stronger and compact reformulation. \citet{zhang-2018-ambig-chanc} derive a submodular relaxation of their DRCCP reformulation for a general binary packing problem and apply extended polymatroid inequalities. \citet{zhang-2020-branch} integrate various algorithmic techniques, including coefficient strengthening and structure-aware reformulation, into a branch-and-price algorithm to solve a bin packing problem.

Tractable reformulations for distributionally robust \emph{joint} chance constraints, i.e., constraint \eqref{eq:dr-ccp-cons} with $m \geq 2$, are much scarcer than for individual chance constraints. Indeed,~\citet[Section 2.3]{hanasusanto-2017-ambig-joint} show that DRCCP becomes NP-hard if the ambiguity set involves any non-homogeneous dispersion measure (e.g., covariance as in $\cF(\mu, \Sigma)$) or any non-conic support (e.g., a hyperrectangle), or if $T(\omega)$ involves any uncertainty (i.e., if $T(\omega) \neq T_0$ for some data matrix $T_0 \in \bbR^{m \times n}$). Nevertheless, tractable reformulations do exist for ambiguity sets different from $\cF(\mu, \Sigma)$ or for chance constraints less general than \eqref{eq:dr-ccp-cons}. For example,~\citet[Theorem 2]{hanasusanto-2017-ambig-joint} characterize an ambiguity set by the mean value, a positively homogeneous dispersion measure (e.g., MAD), and a conic support of $\omega$, and derive a second-order conic reformulation of constraint \eqref{eq:dr-ccp-cons}, in which $T(\omega) = T_0$.~\citet[Theorem 2]{Xie-2018-dr-opf} consider a two-sided variant of \eqref{eq:dr-ccp-cons} with $m = 2$ and $T_1(\omega) = - T_2(\omega)$ and derive a second-order conic reformulation of constraint \eqref{eq:dr-ccp-cons} with regard to ambiguity set $\cF(\mu, \Sigma)$.~\citet{xie-2016-deter-refor} derive exact and tractable reformulations of \eqref{eq:dr-ccp-cons} with regard to multiple types ambiguity sets, e.g., when $\cF(\beta)$ involves linear moment constraints only (i.e., on the mean value of $\omega$) or when $\cF(\beta)$ consists of a single (possibly nonlinear) moment constraint.~\citet{xie-2019-optim-bonfer} consider a subclass of constraints \eqref{eq:dr-ccp-cons} with separable uncertainties across individual inequalities, i.e., each row of $[T(\omega); r(\omega)]$ involves a different set of uncertain parameters and, correspondingly, a different ambiguity set. They show that, if either $T(\omega)$ or $r(\omega)$ involves no uncertainty, then \eqref{eq:dr-ccp-cons} admits an exact and tractable reformulation by applying the Bonferroni approximation (or union bound; see~\citet{bonferroni1936teoria}).

Various \emph{conservative} approximations for distributionally robust joint chance constraints have been proposed.
\citet{chen-2010-from-cvar} propose to approximate the chance constraint in \eqref{eq:dr-ccp-cons} by using CVaR and subsequently approximate the resulting distributionally robust CVaR (DR-CVaR) constraint via a classical inequality of order statistics. These two layers of approximation lead to a set of second-order conic constraints. Later,~\citet{zymler-2011-distr-robus} show that the second-layer approximation can be circumvented by deriving an exact reformulation of the DR-CVaR constraint, yielding a linear matrix inequality approximation of \eqref{eq:dr-ccp-cons}. The approximations of~\cite{chen-2010-from-cvar} and~\cite{zymler-2011-distr-robus} can both be further improved by tuning certain scaling parameters. Unfortunately, it appears to be difficult to simultaneously optimize such scaling parameters and the decision $x$ in DRCCP.~\citet{cheng2014distributionally} obtain a different approximation from that of~\cite{zymler-2011-distr-robus} when different rows of $T(\omega)$ are independent.

\begin{figure}[t]
        \centering
        \begin{subfigure}[b]{0.45\textwidth}
                \includegraphics[width=\textwidth]{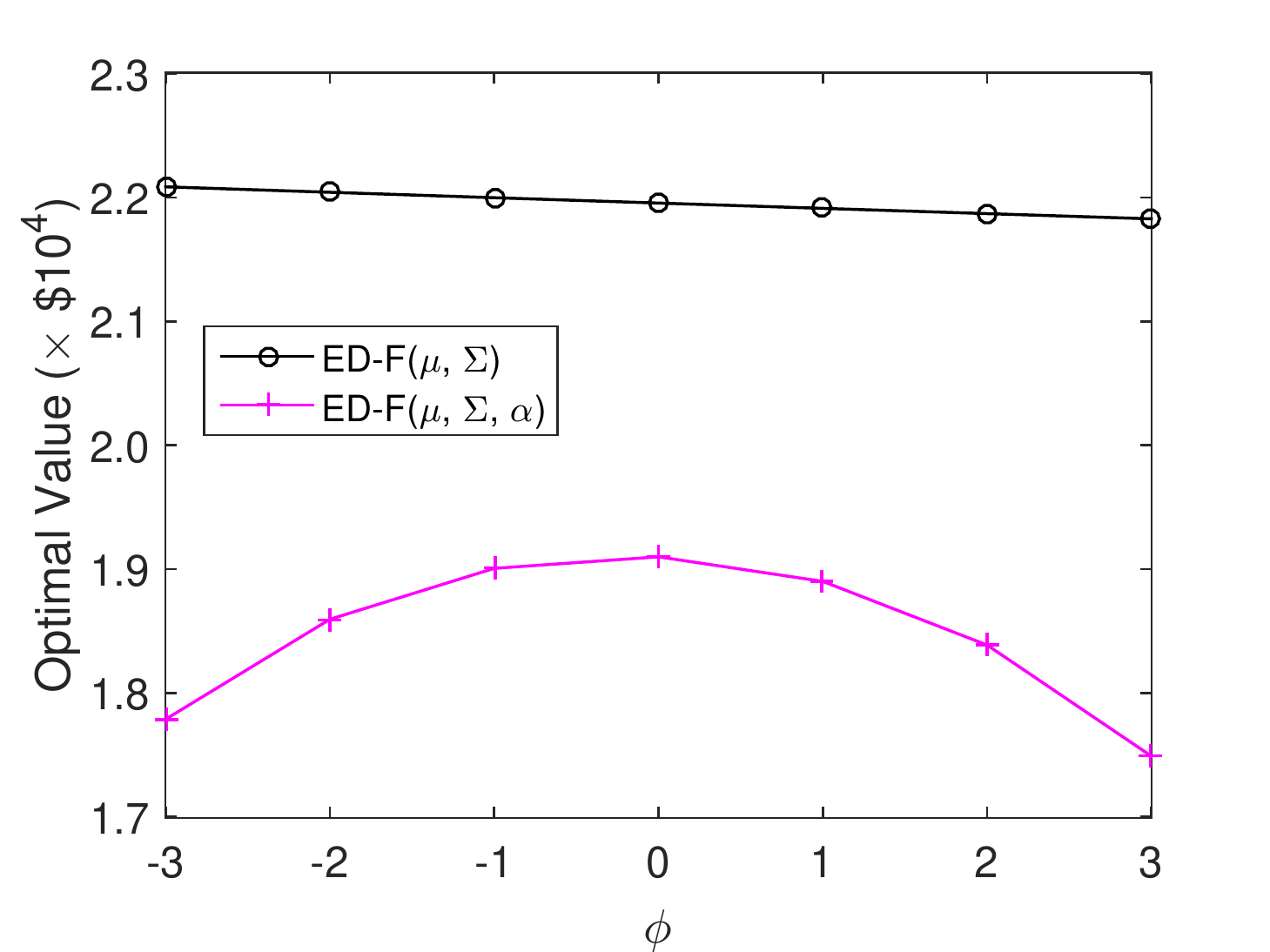}
                \caption{Optimal Value vs. $\phi$}
                \label{fig-ovc-comparison}
        \end{subfigure} \quad %
        \begin{subfigure}[b]{0.45\textwidth}
                \includegraphics[width=\textwidth]{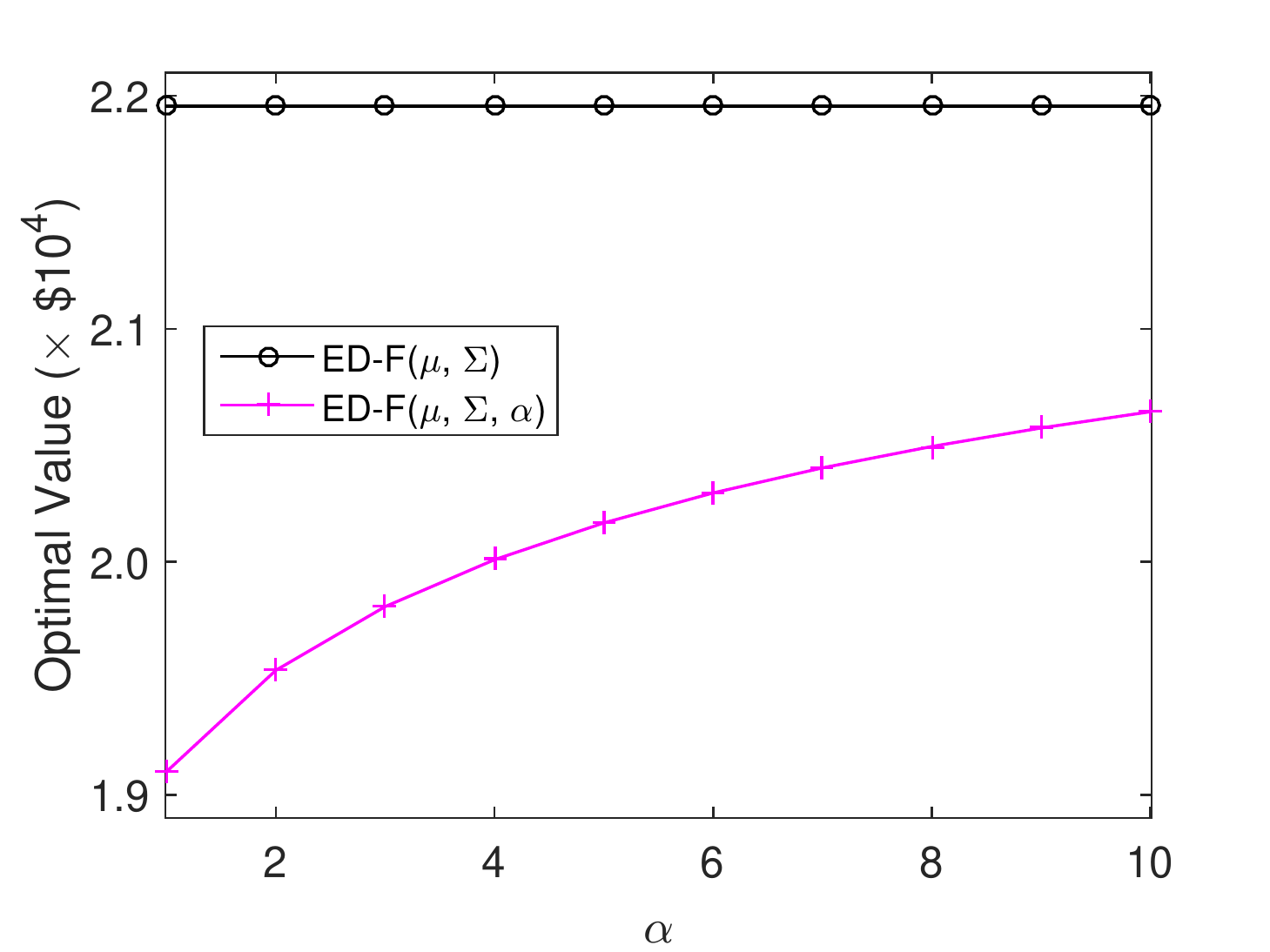}
                \caption{Optimal Value vs. $\alpha$}
                \label{fig-alpha-comparison}
        \end{subfigure}
        \caption{Optimal values of ED-$\cF(\mu, \Sigma)$ and ED-$\cF(\mu, \Sigma, \alpha)$ with various $\phi$ and $\alpha$ {\color{blue}(adapted from Figure 3 of~\cite{li-2019-ambig-risk})}}
\end{figure}

In Figs.~\ref{fig-ovc-comparison}--\ref{fig-alpha-comparison}, we summarize a case study of a distributionally robust chance-constrained economic dispatch (ED) problem that appears in~\citet{li-2019-ambig-risk} to demonstrate the difference between $\cF(\mu, \Sigma)$ and an alternative ambiguity set that incorporates $\alpha$-unimodality into $\cF(\mu, \Sigma)$, denoted by $\cF(\mu, \Sigma, \alpha)$. Their case study uses the IEEE 30-bus system and incorporates two uncertain parameters, representing prediction errors of the forecast power outputs at two wind farms. The formulation and parameters of this problem can be found in~\cite[Section 5.1]{li-2019-ambig-risk}. In particular, we assume that the uncertainties are $\alpha$-unimodal with a mode at $[0,0]^{\top}$ and have a mean value $\mu = \phi[1,1]^{\top}$ with $\phi \in \{-3, -2, \ldots, 3\}$. In Fig.~\ref{fig-ovc-comparison}, we compare the optimal values of ED with regard to $\cF(\mu, \Sigma)$ and that of ED with regard to $\cF(\mu, \Sigma, \alpha)$ with $\alpha = 1$ and various $\phi$ values. From this figure, we observe that the optimal value of ED-$\cF(\mu, \Sigma)$ is consistently higher than that of ED-$\cF(\mu, \Sigma, \alpha)$. This confirms that incorporating unimodality into the ambiguity set makes DRCCP less conservative. In Fig.~\ref{fig-alpha-comparison}, we compare the optimal values of ED-$\cF(\mu, \Sigma)$ and ED-$\cF(\mu, \Sigma, \alpha)$ with $\phi = 0$ and various $\alpha$ values. From this figure, we observe that, although the discrepancy between ED-$\cF(\mu, \Sigma)$ and ED-$\cF(\mu, \Sigma, \alpha)$ declines as $\alpha$ increases, the convergence is sub-linear (in fact, it takes place when $\alpha$ exceeds $10^4$). This demonstrates the significant influence of unimodality upon the ambiguity set and the corresponding DRCCP.

The case study just described highlights the utility of available distribution information in reducing the degree of conservatism. In this regard, moment ambiguity sets are known to be more conservative than their counterparts based on discrepancy measures (e.g., a Wasserstein ambiguity set) when more data samples are available. On the other hand, there is a trade-off between conservatism and tractability---unlike with moment-based ambiguity sets, DRCCP with a Wasserstein ambiguity set is not polynomially solvable in general~\cite{xie-2020-bicrit-approx}. However, there have been recent developments in MIP formulations for DRCCP under Wasserstein ambiguity, which we describe in the next section.

\subsection{Wasserstein ambiguity}\label{sec:wass}

Due to its desirable statistical properties,  the so-called  \emph{Wasserstein} ambiguity set has witnessed an explosion of interest. 
 Wasserstein ambiguity set  $\cF(N,\theta)$ is defined as the $\theta$-radius Wasserstein ball of distributions on $\bbR^d$ around  the empirical distribution $\bbP_N$.  This is defined as
\[ d_W(\bbP,\bbP') := \inf_{\Pi} \left\{ \bbE_{(\omega,\omega') \sim \Pi}[\|\omega - \omega'\|] : \Pi \text{ has marginal distributions } \bbP, \bbP' \right\}, \]
where  the \emph{1-Wasserstein distance}, based on a norm $\|\cdot\|$, between two distributions $\bbP$ and $\bbP'$ is used.
The Wasserstein ambiguity set is then defined as $\cF(\bbP_N,\theta) := \left\{ \bbP : d_W(\bbP_N,\bbP) \leq \theta\right\}.$
Given a decision $x \in \cX$ and random realization $\omega \in \bbR^d$, we first define a safety set, $\cS(x)$, of outcomes such that $\cS(x)=\{\omega\in \Omega: x\in \cP(w)\}$.   The  distance from $\omega$ to the unsafe set is 
\begin{equation}\label{eq:distance}
\dist(\omega,\cS(x)) := \inf_{\omega' \in \bbR^d} \left\{ \|\omega - \omega'\| : \omega' \not\in \cS(x) \right\}.
\end{equation}

\citet[Theorem 3]{chen-2018-data-driven} and \citet[Proposition 1]{xie-2019-distr-robus} show that the formulation 
\begin{subequations}\label{eq:cc-distance-formulation}
\begin{align}
\min_{x,v,u} \quad & c^\top x\notag\\ 
&\quad   x \in \cX, \ v \geq 0, u_i \geq 0, \ i \in [N], \\
&\quad \dist(\omega_i,\cS(x)) \geq v - u_i, \ i \in [N], \label{eq:drcc-nonconvex}\\
&\quad \epsilon\, v \geq \theta + \frac{1}{N} \sum_{i \in [N]} u_i
\end{align}
\end{subequations}
is an equivalent formulation of \eqref{eq:dr-ccp}, by using the dual representation for the worst-case probability $\bbP[x \not\in \cP(\omega)]$ under the Wasserstein ambiguity set $\bbP \in \cF(\bbP_N,\theta)$ provided in \cite{BlanchetMurthy2019,gao2016distributionally,MohajerinEsfahaniKuhn2018}. (See also \citet{HotaEtAl2019} for a deterministic non-convex reformulation of \eqref{eq:dr-ccp} and CVaR-based inner approximation of \eqref{eq:dr-ccp} for certain safety sets.)

Note that formulation \eqref{eq:dr-ccp}  is non-convex due to  constraint \eqref{eq:drcc-nonconvex}. However, for certain safety sets $\cS(\cdot)$, MIP reformulations are possible  \cite{chen-2018-data-driven,xie-2019-distr-robus,ji-2018-data-driven}. Therefore, we can once again formulate a deterministic equivalent model to be able to solve it using off-the-shelf optimization software.

\subsubsection{RHS Uncertainty}

In this section, we consider  joint chance constraints with RHS uncertainty under certain common form of a safety set. In particular, let 
\begin{equation}\label{eq:joint-safety-set}
    \cS(x) := \left\{ \omega :~  Tx\ge r(\omega) \right\},  \\
\end{equation}
where $r(\omega):=B\omega+e$, for a given  an $m\times d$ data matrix $B$, $e\in \bbR^m$, and $T$ is a given $m\times n$ data matrix. For $m=1$ (resp.\ $m>1$), we say that the problem is an individual (resp.\ joint) chance-constrained problem with RHS uncertainty. Let $T_j$ and $B_j$ be a row vector of appropriate dimension corresponding to the $j$th row of $T$ and $B$, respectively.  In this case, the distance function is evaluated as \cite{chen-2018-data-driven}
\begin{equation}
\dist(\omega,\cS(x))  = \max\left\{ 0,\ \min_{j \in [m]} \frac{T_j x - B_j \omega - e_j }{\|B_j\|_*} \right\},
\end{equation}
where $\|\cdot\|_*$  is the dual norm. We can then introduce binary variables, $z$, to capture the non-convex constraint \eqref{eq:drcc-nonconvex} to arrive at the mixed-integer \emph{linear} program \cite[Proposition 2]{chen-2018-data-driven}
\begin{subequations}\label{eq:joint}
\begin{align}
\min\limits_{z, u, v, x} \quad & c^\top x\label{joint:obj}\\
\text{s.t.}\quad & z \in \{0,1\}^N,\ v \geq 0, \ u_i \geq 0, \ i\in [N],\ x \in \mathcal{X},\label{joint:vars}\\
& \epsilon\, v \geq \theta + \frac{1}{N} \sum_{i \in [N]} u_i,\label{joint:conic}\\
& M (1-z_i) \geq v-u_i, \quad i \in [N],\label{joint:bigM1}\\
& \frac{T_j x - B_j \omega_i - e_j }{\|B_j\|_*}+ M_i z_i\geq v-u_i, \quad i \in [N],\ j \in [m],\label{joint:bigM2}
\end{align}
\end{subequations}
where $M_i, i \in [N]$ is a sufficiently large Big-M coefficient. 

A few remarks are in order. The computational studies  of \cite{chen-2018-data-driven,xie-2019-distr-robus} indicate that this MIP reformulation is  difficult to solve in certain cases---state-of-the-art solvers  terminate with large optimality gaps  after an hour time limit. To address this challenge, \citet{ho-nguyen-2020-distr-robus}  propose a number of results that make an order of magnitude improvement in the solution times. Note that formulation \eqref{eq:joint} is not immediately amenable to the improvements we described for the SAA counterpart. For example, constraints \eqref{joint:bigM2} do not have the mixing structure that the SAA counterpart benefited greatly from. In particular, the continuous variables $u_i$ are not shared across scenarios, whereas the mixing set requires common continuous variables. On the other hand, as argued in \cite{ho-nguyen-2020-distr-robus}, the SAA counterpart is a relaxation of \eqref{eq:joint}. By making a key observation that relates the nominal SAA problem for $\bbP_N$ to formulation \eqref{eq:joint}, \citet{ho-nguyen-2020-distr-robus} give a stronger formulation and valid inequalities based on the same set of binary variables $z$. Furthermore, this strengthening does have the mixing structure.  They also use pre-processing techniques to reduce the formulation size drastically. 
On a  related note,  \citet{ji-2018-data-driven} give a different MIP formulation of \eqref{eq:dr-ccp} under Wasserstein ambiguity under additional assumptions on the support of $\omega$.

\subsubsection{LHS uncertainty}

In this section, we consider  joint chance constraints with RHS uncertainty under certain common form of a safety set. In particular, let 
\begin{equation}\label{eq:joint-safety-set}
    \cS(x) := \left\{ \omega :~  T(\omega)x\ge r(\omega) \right\},  \\
\end{equation}
where $r_j(\omega):=b^\top \omega^j +e_j, j\in [m]$, for a given vector $b\in \bbR^\kappa$, $\omega^j, j\in[m]$ is a projection of $\omega$ to a $\kappa$-dimensional vector, and $e\in \bbR^m$. Also, let the $j$th row of $T(\omega)$ be given by $T_j(\omega):=\omega^\top A+T_j$ for some $n\times \kappa$ data matrix $A^\top$ and $T\in \bbR^{m\times n}$.  In this case, the distance function is measured by
\begin{equation}
\dist(\omega,\cS(\vx)) = \max\left\{ 0,\ \min_{j \in [m]} \frac{T_j(\omega)x-r_j(\omega)}{\|A^\top x -b \|_*} \right\},\label{eq:distance-linear}
\end{equation}
We can then introduce binary variables, $z$ to represent the non-convex constraint \eqref{eq:drcc-nonconvex} and make a transformation of variables to arrive at the mixed-integer \emph{conic} program (\cite[Theorem 2]{xie-2019-distr-robus} and \cite[Proposition 1  (for $m=1$)]{chen-2019-distr-robus}
\begin{align*}
\min\limits_{z, u, v, x} \quad & c^\top x\\
\text{s.t.}\quad & z \in \{0,1\}^N,\ v \geq 0, \ u_i \geq 0, \ i\in [N],\ x \in \mathcal{X},\\
& \epsilon\, v \geq \theta \|A^\top x -b \|_* + \frac{1}{N} \sum_{i \in [N]} u_i,\\
& M_i (1-z_i) \geq v-u_i, \quad i \in [N],\\
& T_j(\omega_i)x-r_j(\omega_i) + M_i z_i\geq v-u_i, \quad i \in [N],\ j \in [m],
\end{align*}
where $M_i, i \in [N]$ is a sufficiently large Big-M coefficient, under the assumption that $ A^\top x\neq b $ for any $x \in \cX$. This assumption can be relaxed with appropriate safeguards as described in \cite{chen-2019-distr-robus,xie-2019-distr-robus,ho-nguyen-2020-stron-formul}.

As in the case of SAA, the computational studies show that the LHS uncertainty case is a more challenging case than the RHS uncertainty.  First, the resulting formulation is no longer linear, but conic. Furthermore, the coefficients of the common variables $x$ are scenario-dependent  unlike the RHS uncertainty case. So it is not clear if similar enhancements that \citet{ho-nguyen-2020-distr-robus} performed for the RHS uncertainty case can be done here. To this end, \citet{ho-nguyen-2020-stron-formul}  establish the link between the DRCCP and its SAA counterpart for the LHS case to identify \revise{valid inequalities based on mixing sets} and strengthen the formulation. This results in significant improvements in the performance of the resulting MIP formulation. Distributionally robust variants of the resource planning problem (described in Section~\ref{sec:two-stage-ccp})  with $N=100$ that are unsolvable or terminate with high end gaps (40-80\%) with the original formulation are now solvable or have much small end gaps (<15\%) with the enhancements proposed in \cite{ho-nguyen-2020-stron-formul}.

For \emph{combinatorial DRCCPs}, for which the  decision variables
are pure binary,  further strengthening is possible.~\citet{xie-2019-distr-robus} observe the submodularity of the norm and the terms in the distance operator, and propose the use of polymatroid inequalities to strengthen the formulation. They report significant improvements in the performance of the resulting algorithm. \citet{Kilinc-Karzan2020-conicbinary} show how the polymatroid inequalities derived from the conic constraint can be generalized to the case of mixed-binary decisions. In addition, \citet{shen-2020-cc-set-cover-wasserstein} derive polymatroid inequalities when the random parameters are binary-valued and show how these inequalities can be further strengthened via mixing and lifting schemes. In a related line of work,~\citet{wang-2019-solut-approac} consider an assignment problem and derive lifted cover inequalities based on a bilinear reformulation of their DRCCP.

\emph{Conservative approximations} for  DRCCP with Wasserstein ambiguity are related to their SAA counterparts described in Section \ref{sec:approx}.
The approach of \citet{erdogan-2005-ambig-chanc} may be seen as a (robust) scenario approximation counterpart of \cite{CampiGaratti2008,calafiore-2004-uncer-convex-progr} 
with similar sample complexity results when the uncertainty set is defined by a Prohorov metric, which is related to a Wasserstein metric. 
Furthermore, for distributionally robust CCPs under Wasserstein ambiguity \cite{HotaEtAl2019} give an approximation based on a CVaR interpretation of the reformulation \cite[see, also,][for this and two other approximations based on the scenario approximation and VaR approximation]{xie-2019-distr-robus}.

\section{Concluding Remarks}\label{sec:conclude}

In this survey, we reviewed \revise{reformulations of CCPs based on sampling and distributional robustness, when there is limited distributional information}. We described the trade-offs between tractability and  conservatism of the corresponding optimization models, as well as the trade-offs between the amount of distributional information used  and over-fitting. There is some theoretical guidance on selecting sample sizes or other design parameters, such as the Wasserstein ball radius. However, this guidance is conservative, and instead the parameter choices are made and statistically verified using out-of-sample tests and cross-validation, in practice. 
There are many opportunities that arise from the recent developments in CCP models.  As we outlined, these models often lead to mixed-integer \revise{(conic)} formulations, which optimization software is now able to handle in modest sizes.   \revise{Coupling the novel mixed-integer conic CCP models   with parallel developments in strengthening mixed-integer conic formulations \cite{atamturk2010conic-mir,Kilinc-Karzan2020-conicbinary,atamturk2019conic-quad, xie-2019-distr-robus,zhang-2018-ambig-chanc,atamturk2008polymatroids} will  likely } enable the solution of large-scale problems before resorting to conservative approximations. Such strengthening approaches often exploit hidden submodularity---a recurring structure in many reformulations we discussed.
 Approximations continue to play an important role in applications where faster solution times are needed. In such cases, it is of interest to be able to provide some performance guarantees. In this regard, recent research in deriving strong relaxations and approximation algorithms for structured problems is promising. 

\revise{We have focused on CCPs where the distribution is unknown in closed form, but either a sample or other limited  information (e.g., moments, unimodality) is available. We refer the reader to \cite{HS08,vAHM10,CHENG2012325,CHENG2013597, Cheng2015, vA15, vAdO16, VANACKOOIJ201199, van2019eventual} and references therein for CCPs under continuous distributions with more knowledge of the distributions (e.g., Gaussian, known copul\ae),  and~\cite{shen2021convex} for DRCCPs with Wasserstein ambiguity and a log-concave reference distribution (in place of an empirical distribution). 
By exploiting  the known properties of the continuous distribution, if available, 
convexity of the feasible region may be established, tractable reformulations and specialized algorithms such as generalized Benders decomposition can be developed.}

We have primarily discussed single- or two-stage problems in this survey. Conceptually, one can also envision CCPs with multiple decision epochs. \citet{Minjiao} consider multi-stage CCPs and give valid inequalities for the SAA reformulation. \citet{LS04} consider a multi-stage  problem under a finite discrete demand distribution, and propose a model wherein non-anticipativity is enforced only for the scenarios that meet the desired service constraint. The authors propose a branch-and-price algorithm, for the resulting formulation. \citet{ANDRIEU2010,Gonzales-Grandon2020}, and references therein, consider problems with dynamic chance constraints, and propose solution methods under certain continuous distributions.  
 \citet{Merakli2019} consider the risk associated with parameter uncertainty in infinite-horizon Markov decision processes, and formulate this problem using a chance-constrained optimization framework.  Models and methods for multi-stage CCPs are sparser due to their inherent difficulty  not only in modeling, by taking into account the time consistency of solutions,  but also in designing scalable solution methods. This is an area of further research. 

\revise{Throughout the survey, we highlighted major trade-offs in comparing alternative methods: tractability, conservatism of solutions from approximations, and ease of implementation. The particular choice of a solution method to choose will depend on the decision-maker's risk preferences, available computational budget, and available resources for implementation of special-purpose methods.}

In closing, we believe that the developments in easy-to-implement reformulations will usher in new and exciting applications of CCPs, given the increasingly uncertain conditions of operations in various sectors (extreme weather, autonomous devices, renewable power, pandemics, political unrest, etc.).

\section*{Acknowledgments}
\revise{We thank the editors and two anonymous referees whose comments improved the positioning of this survey within the vast area of CCPs. We also thank Weijun Xie for his comments on an earlier version.} 
Simge K\"u\c{c}\"ukyavuz is supported, in part, by ONR grant N00014-19-1-2321 and NSF grant 2007814. Ruiwei Jiang is supported, in part, by NSF grant ECCS-1845980.

\appendix

\section{Applications} \label{sec:app}


In this section, we review a few recent and active applications of CCPs in practice---this is not meant to be an exhaustive list. 

\paragraph{Finance.} Chance constraints (or equivalently, VaR as defined in \eqref{vardef}) have been applied in finance to control risks. \citet{linsmeier2000value} provide motivation of using VaR as a risk measure in significant volatile financial markets. VaR has been widely adopted (e.g., by the US Securities and Exchange Commission) as a method of quantifying risks. \citet{lemus1999portfolio}, \citet{ghaoui-2003-worst-case}, \citet{natarajan-2008-var}, \citet{zymler2013worst}, \citet{HUANG2014243}, \citet{yao2015smooth}, \citet{ccetinkaya2015data}, \citet{BARRIEU2015546}, \citet{LOTFI2018556}, \citet{li2018worst-range-var}, and~\citet{ji2018risk} apply VaR and worst-case VaR (analogous to the distributionally robust chance constraints) in finance via mathematical optimization. In addition, \citet{napat-2016-growth} and \citet{choi2016multi} apply chance constraints in multi-period portfolio optimization.

\paragraph{Healthcare.} Chance constraints find applications in appointment scheduling (e.g.,~\citet{deng2016decomposition}), surgery planning (e.g.,~\citet{deng2019chance}, \citet{wang2017distributionally}, and~\citet{zhang2018solving}), operating room planning (e.g.,~\citet{wang2019chance}, \citet{wang-2019-solut-approac}, and~\citet{najjarbashi2020decomposition}), vaccine allocation (e.g.,~\citet{Tanner2010}), and social distancing during a pandemic (e.g.,~\citet{duque2020timing}), among others.

\paragraph{Power Systems.} \citet{chang-li-2011-chance-opf}, \citet{bienstock2014chance}, \citet{zhang-shen-2017-drcc-opf}, \citet{duan-2018-drcc-opf-wasserstein}, \citet{lubin-2016-chance,Lubin2019} \citet{dallanese-2017-chance-ac-opf}, \citet{Xie-2018-dr-opf}, \citet{li2018distributionally}, and~\citet{li2019distributionally} study chance-constrained variants of the optimal power flow problem. \citet{ozturk-2004-chance-uc}, \citet{pozo-2013-chance-uc-n-k}, and \citet{WGW12} consider chance constraints in the unit commitment problem. \citet{vrakopoulou-2013-probabilistic}, \citet{pozo-2013-chance-uc-n-k}, and~\citet{wu-2014-chance} apply chance constraints to schedule electricity systems in face of random outages and contingencies. \citet{liu-2011-optimal-dg}, \citet{liu-2018-ev}, \citet{ravichandran-2018-chance-microgrid}, and~\citet{zhang2020values} employ chance constraints to model an integrated system of power grid and electric vehicles. Other power system applications include coordinated load control (e.g.,~\citet{zhang-shen-2017-drcc-opf} and~\citet{zhang2019distributionally}), power grid topology control (e.g.,~\citet{qiu-2015-chance-switching} and~\citet{mazadi-2009-chance-expansion}), and hydro power plant scheduling (e.g.,~\citet{wu-2008-chance-hydro} and~\citet{Lodi2019}). We refer the reader to a recent survey \cite{vanAckooij-2011-survey} and references therein for a more detailed review of CCP in energy management. 

\paragraph{Transportation and Routing.} \citet{dinh2018exact}, \citet{moser-2018-flexible}, \citet{pelletier2019electric}, \citet{du-2020-cooperative}, \citet{wu-2020-safe}, \citet{muraleedharan-2020-chance}, \citet{ghosal-2020-drcc-vrp}, and~\citet{florio-2021-chance-vrp} study chance constraints in the optimal route design for vehicles (also see~\citet{cordeau2007vehicle}). \citet{balckmore-2011-chance}, \citet{farrokhsiar-2011-unscented}, \citet{banerjee-2011-regression}, \citet{du-toit-2011-chance}, \citet{arantes-2019-collision-chance}, \citet{castillo-lopez-2020-chance-obstacle}, and~\citet{oh-2020-chance-path} study chance constraints to find paths for robots while avoiding obstacles.

\paragraph{Supply Chain, Logistics, and Scheduling.} \citet{wang2007beta}, \citet{song-2013-branc-and}, \citet{Hong2015}, \citet{Elci2018}, \citet{ElciNB2018}, and~\citet{Noyan2019} employ chance constraints in the design of networks for logistics and humanitarian relief. \citet{LR07}, \citet{Murr2000}, \citet{Minjiao}, and~\citet{LK18} apply chance constraints in logistics. \citet{GLT10} study chance constraints in the staffing of call centers. \citet{cohen2019overcommitment} apply chance constraints to cloud computing. \citet{lu2020non} apply chance constraints in non-profit resource allocation.

\paragraph{Wireless Communication.} \citet{li-2010-slow-adap},~\citet{soltani-2013-chance},~\citet{mokari-2016-robust}, and~\citet{xu-2016-energy-chance} apply chance-constrained programming to accommodate the data rate requirement in orthogonal frequency division multiple access (OFDMA) systems. \citet{ma-2013-chance-beamforming} and~\citet{li2014distributionally} apply chance constraints on the beamforming problem in communication networks.

\section{Other Surveys on CCPs}\label{sec:surveys}
\revise{
In Table \ref{tab:survey}, we provide a list of other survey papers and monographs on CCPs, their foci, and the sections in which we referred to them in this paper.

\begin{table}[thb]
\caption{Reviews on CCPs} \label{tab:survey} \centering
\resizebox{\textwidth}{!}{
\begin{tabular}{c|c|c}
Reference&Focus&Section\\\hline
\cite{BL97,KW94,P95,prekopa2003ccp} & Monographs on work prior to 2000  & \ref{sec:intro}\\
\cite{D09.ch4,Dentcheva06} & Monograph, special distributions, $(1-\epsilon)$-efficient points  & \ref{sec:intro},~\ref{sec:saa-alt-form} \\
\cite{NEMIROVSKI2012}& Approximations of individual chance constraints & \ref{sec:approx} \\
\cite{vanAckooij-2011-survey} & Special distributions (e.g., i.i.d. Gaussian); energy applications  & \ref{sec:def}, Appendix~\ref{sec:app}\\
\cite{Ahmed2018} & Relaxations and approximations of CCP & \ref{sec:approx}\\
\cite{Lejeune2018} & Relaxations and approximations of joint CCPs & \ref{sec:approx}\\
\cite{AS08} & SAA approaches, RHS uncertainty & \ref{sec:saa-rhs}\\
This paper & SAA and DRO-based reformulations, two-stage CCPs, applications & \ref{sec:saa-rhs},~\ref{sec:saa-joint},~\ref{sec:two-stage-ccp},~\ref{sec:dro}, Appendix~\ref{sec:app}
\end{tabular}
}
\end{table}
}

\section{Preliminaries}\label{sec:prelim}
In this section we review some preliminaries.

\subsection{Value-at-risk and Conditional value-at-risk}
Here we  present two relevant
definitions pertaining to the risk associated with a univariate random variable that are used in our discussion.   We refer the reader to  
\cite{Pflug00,Pflug07,Rockafellar07} for a more detailed treatment of these risk measures. 

\begin{definition} \label{def:var}
For a univariate random variable $X$, with cumulative distribution function $F_X$, 
the \emph{value-at-risk} ($\var$)
at confidence level $(1-\epsilon)$, also known as $(1-\epsilon)$-quantile, is  given by:
\begin{equation}\label{vardef}\var_{1-\epsilon}(X)=\min\{\eta~:~ F_X(\eta)\geq
1-\epsilon\}.\end{equation}
\end{definition} 
It follows from \eqref{vardef} that, for any $x \in \bbR$, the inequalities $\var_{1-\epsilon}(X) \leq x$ and $\bbP(X \leq x) \geq 1 - \epsilon$ are equivalent. That is, a chance constraint on random variable $X$ can be equivalently represented as a constraint on its $\var$.
\begin{definition}[\cite{Rockafellar00,Rockafellar02}]\label{def-cvar}
The \emph{conditional value-at-risk} (CVaR)
at confidence level $(1-\epsilon)\in(0,1]$ is given by
\begin{equation}\label{cvardef}\cvar_{1-\epsilon}(X)=\min\left\{\eta+\frac{1}{\epsilon}\E\left([X-\eta]_+\right)
~:~\eta\in\R\right\},\end{equation}
where $(a)_+:=\max\{0,a\}$.
\end{definition}
It is well known that
the minimum in definition \eqref{cvardef} is attained at the
$\var$
at confidence level $(1-\epsilon)$. 
 CVaR, introduced by
\citet{Rockafellar00}, satisfies the axioms of coherent risk measures, such as law invariance and
sub-additivity, as defined in \cite{Artzner99}. It has other desirable properties, such as tractability---for finite distributions, CVaR can be formulated as a linear program and embedded in an optimization model \cite{Rockafellar07}. 
More precisely, suppose $X$ is a random variable with 
realizations $X_1,\ldots,X_N$ and corresponding probabilities
$p_1,\ldots,p_N$. \ignore{Throughout, for $a\in\bbZ_+$, let $[a]:=\{1,\dots,a\}$.} 
The optimization problem in \eqref{cvardef} can equivalently be
formulated as the linear program (LP):
\begin{equation}\label{cvar_classic}\min \left\{\eta+\frac{1}{\epsilon}\sum_{i\in [N]}p_i
w_i~:~w_i\geq X_i- \eta,~\forall~i\in [N],\quad w\in
\R_+^{N}\right\}.\end{equation}

Furthermore,  let $\rho$ denote an ordering of the realizations such that $X_{\rho_1}\leq X_{\rho_2}\leq \cdots \leq X_{\rho_N}$. Then, for a given confidence
level $\epsilon\in (0,1]$ we have
\begin{equation}
\var_{1-\epsilon}(X)=X_{\rho_q}, \text{ where } q=\min\left\{j \in
[N]:\sum_{i\in [j]}p_{\rho_i}\geq
1-\epsilon\right\}.\label{varknapsack}\end{equation}

\revise{
\subsection{Submodularity and polymatroid inequalities} \label{sec:submod}
Here we review the concept of submodularity and valid inequalities for the epigraph of a submodular function. 
\begin{definition}
A set function $g: 2^N\to \mathbb{R}$ is submodular if 
		\[
		g(A)+g(B)\geq g(A\cap B) + g(A\cup B)\quad \forall A,B\subseteq[N].
		\]
\end{definition}

 Given a submodular (set) function $g:2^{N}\to\mathbb{R}$, the \emph{extended polymatroid} of $g$ is
	\begin{equation*}
	EP_g:=\left\{\pi\in\mathbb{R}^N:~ \sum_{j\in V}\pi_j\leq g(V),~\forall V\subseteq[N]\right\}.
	\end{equation*}
	
 By slightly abusing notation, for any $z\in \{0,1\}^N$, $g(z)$ is equivalent to $g(S)$, where $S$ is the support of $z$, i.e., for all $i \in [N]$, $z_i = 1$ if and only if $i \in S$. Given a submodular function $g:\{0,1\}^n\to\mathbb{R}$, consider its epigraph
	\begin{equation*}
	Q_g:=\left\{(t,z)\in\mathbb{R}\times\{0,1\}^N:~t\geq g(z)\right\}.
	\end{equation*}

	\begin{theorem}\cite{lovasz1983submodular, atamturk2008polymatroids} \label{thm:polymat}
	The convex hull of $Q_g$ is given by
\begin{equation*}\left\{(t,z)\in\mathbb{R}\times[0,1]^N:~{t\geq\pi^\top  z+g(\emptyset)},~\forall \pi\in EP_{g-g(\emptyset)}\right\}. \label{eq:epi}
\end{equation*}
	\end{theorem}

 The inequalities ${t\geq \pi^\top z+g(\emptyset)}$ for $\pi\in EP_{g-g(\emptyset)}$ are referred to as the \emph{polymatroid inequalities} of $g$.
	
 Polymatroid inequalities can be separated  in $O(N\log N)$ time.

}

\revise{
\section{Additional Valid Inequalities for Joint CCPs with RHS Uncertainty} \label{sec:blending}

Consider the SAA formulation of a joint chance constraint with RHS uncertainty:
\begin{align}
& t_j \geq r_{i,j}(1 - z_i), \quad \forall i \in [N], \ \forall j \in [m], \label{eq:intmix-1} \\
& \sum_{i \in [N]} p_i z_i \leq \epsilon.  \label{eq:intmix-2}
\end{align}
One can follow the description in Section~\ref{sec:saa-rhs} to obtain valid inequalities for the set
\[
\cM_j = \left\{(t,z)\in\R_+\times\{0,1\}^N:~ t_j\geq r_{i,j} (1-z_i),~\forall i\in[N], \sum_{i\in[N]} p_i z_i \leq \epsilon\right\},
\]
that is, the mixing set with regard to an individual inequality for each $j \in [m]$. Furthermore,~\citet{zhao2017joint-knapsack} propose a procedure to \emph{blend} these mixing inequalities into a stronger one for the intersection of all $\cM_j$'s (denoted as $\bigcap\limits_{j\in[m]}\cM_j$, with a slight abuse of notation) given by the system \eqref{eq:intmix-1}--\eqref{eq:intmix-2}. To describe this procedure concretely, we define the following notation.
\begin{definition}{(Adapted from Definition 1 in~\citet{zhao2017joint-knapsack})}
First, for all $j \in [m]$, let $\{\langle i \rangle_j: i \in [N]\}$ denote a permutation of the set $[N]$ such that $r_{\langle 1 \rangle_j, j} \geq r_{\langle 2 \rangle_j, j} \geq \cdots \geq r_{\langle N \rangle_j, j}$. Accordingly, we define $\nu_j$ such that
\[
\sum_{i=1}^{\nu_j} p_{\langle i \rangle_j} \leq \epsilon \quad \text{and} \quad \sum_{i=1}^{\nu_j + 1} p_{\langle i \rangle_j} > \epsilon.
\]
Second, as briefly described in Section \ref{sec:saa-rhs}, we let $\pi$ denote another permutation of $[N]$ such that $p_{\pi_1} \leq p_{\pi_w} \leq \cdots \leq p_{\pi_N}$ and define $\varphi$ such that
\[
\sum_{i=1}^{\varphi} p_{\pi_i} \leq \epsilon \quad \text{and} \quad \sum_{i=1}^{\varphi+1} p_{\pi_i} > \epsilon.
\]
In other words, $\varphi$ provides a valid cardinality bound on the binary variables. 
Notice that $\varphi$ does not depend on the index $j$ and $\varphi \geq \nu_j$ for all $j \in [m]$. Third, for $\theta \in [N]$ and, for each $j \in [m]$, $m_j \in [\nu_j]$ and $q_j \in \{0\}\cup[\varphi-m_j]$, we define
\begin{itemize}
\item a set $U_j = \{u_{j1}, u_{j2}, \ldots, u_{ja_j}\} \subseteq [m_j]$ with $u_{j1} < u_{j2} < \cdots < u_{ja_j}$;
\item a set $L_j$ with $|L_j| = q_j$ and $\theta \in L_j$;
\item a sequence of integers $\{s_{jk}: k \in [q_j + 1]\} \subseteq \{0\}\cup[\nu_j - m_j + 1]$ such that
\begin{itemize}
\item $0 \leq s_{j1} \leq \cdots \leq s_{j,q_j+1} = \nu_j-m_j+1$;
\item $L_j \subseteq \{m_j + s_{j1} + 1, \ldots, N\}$;
\item there exists a permutation $\{\ell_{j1}, \ell_{j2}, \ldots, \ell_{j,q_j}\}$ of $L_j$ with $\ell_{j1} = \theta$ and $\ell_{jk} \geq m_j + \min\{s_{jk}+1, s_{j,k+1}\}$ for all $k \in [q_j]$;
\item there exists a permutation $\{h_{j1}, h_{j2}, \ldots, h_{j,q_j}\}$ of $L_j$ such that $p_{\langle h_{j1} \rangle_j} \geq \cdots \geq p_{\langle h_{j,q_j} \rangle_j}$ and
\[
\sum_{i=1}^{m_j + s_{jk}} p_{\langle i \rangle_j} + \sum_{i=k}^{q_j} p_{\langle h_{ji} \rangle_j} > \epsilon, \quad \forall k \in [q_j].
\]
\end{itemize}
\end{itemize}
\end{definition}
Then, a \emph{strengthened} mixing inequality valid for the system \eqref{eq:intmix-1}--\eqref{eq:intmix-2} is given by Theorem 2.1 of~\citet{zhao2017joint-knapsack} (also see Theorem 6 of~\citet{kuecuekyavuz-2012-mixin-sets}):
\begin{equation}
t_j + \sum_{k=1}^{a_j} \left(r_{\langle u_{jk} \rangle_j,j} - r_{\langle u_{j,k+1} \rangle_j,j}\right)z_{\langle u_{jk} \rangle_j} + \sum_{k=1}^{q_j} \delta_{jk} (1 - z_{\langle \ell_{jk} \rangle_j}) \geq r_{\langle u_{j1} \rangle_j,j}, \label{eq:mixing-individual}
\end{equation}
where
\begin{equation*}
\delta_{jk} = \left\{\begin{array}{ll}
    r_{\langle m_j + s_{j1} \rangle_j,j} - r_{\langle m_j + s_{j2} \rangle_j,j} & \text{if $k=1$} \\[1em]
    \max\left\{\delta_{j,k-1},~r_{\langle m_j + s_{j1} \rangle_j,j} - r_{\langle m_j + s_{j,k+1} \rangle_j,j} - \displaystyle\sum_{i<k:~\ell_{ji} \geq m_j + \min\{1+s_{jk},s_{j,k+1}\}}\delta_{ji}\right\} & \text{for all $2 \leq k \leq q_j$}.
\end{array}\right.
\end{equation*}
Now, suppose that the sets $\langle U_j \rangle_j \setminus \{\theta\}$ for all $j \in [m]$ are mutually disjoint, where $\langle X \rangle_j = \{\langle i \rangle_j:~i \in X\}$ for any set $X \subseteq [N]$, and that $\sum_{i \in \cup_{j \in [m]}\langle [m_j + s_{j1}] \rangle_j}p_i > \epsilon$. Then,   inequalities~\eqref{eq:mixing-individual} can be blended into the following (stronger) valid inequality for the system \eqref{eq:intmix-1}--\eqref{eq:intmix-2} (see Theorem 4.1 in~\citet{zhao2017joint-knapsack}):
\begin{equation}
\sum_{j \in [m]} \frac{1}{\delta_{j1}} \left(t_j + \sum_{k=1}^{a_j} \left(r_{\langle u_{jk} \rangle_j,j} - r_{\langle u_{j,k+1} \rangle_j,j}\right)z_{\langle u_{jk} \rangle_j} + \sum_{k=1}^{q_j} \delta_{jk} (1 - z_{\langle \ell_{jk} \rangle_j}) - r_{\langle u_{j1} \rangle_j,j}\right) \ \geq \ 1 - z_{\theta}. \label{eq:mixing-blended}
\end{equation}

We illustrate the blended inequality~\eqref{eq:mixing-blended} on our numerical example (Example~\ref{ex:sen-ex}; also see Example 3 of~\citet{zhao2017joint-knapsack} and Example 2 of~\citet{kuecuekyavuz-2012-mixin-sets}), in which $m=2$, $N=9$, and $\epsilon=0.4$. By definition, we obtain $k = 6$, $\nu_1 = 3$, and $\nu_2 = 5$. In addition, we obtain a mixing inequality
\[
t_1 + 0.25z_1 + 0.25(1 - z_9) \geq 0.75
\]
for $\cM_1$ with $m_1 = 1$, $q_1 = 1$, $U_1 = \{1\}$, $L_1 = \{9\}$, and $\{s_{11}, s_{12}\} = \{1, 3\}$, and another mixing inequality
\[
t_2 + 0.5z_7 + 0.25(1 - z_9) \geq 2
\]
for $\cM_2$ with $m_2 = 1$, $q_2 = 1$, $U_2 = \{1\}$, $L_2 = \{9\}$, and $\{s_{21}, s_{22}\} = \{4, 5\}$. Then, for $\theta = 9$, we blend these two mixing inequalities to obtain
\[
t_1 + t_2 + 0.25 z_1 + 0.5 z_7 + 0.25(1-z_9) \geq 2.75,
\]
which is valid for $\cM_1 \cap \cM_2$.

There is no known exact separation algorithm for inequalities \eqref{eq:mixing-individual} and \eqref{eq:mixing-blended}. However, effective heuristics are proposed in  \cite{zhao2017joint-knapsack}.

In another line of work, \citet{Kilinc-Karzan2019joint-submod} consider a set of the form
\begin{align*}
& t_j \geq r_{i,j}(1 - z_i), \quad \forall i \in [N], \ \forall j \in [m], \\
& \sum_{j\in [m]} t_j \geq \varepsilon, 
\end{align*}
for a given $\varepsilon$ that provides a lower bound on $\sum_{j\in [m]} t_j$ (e.g., that obtained from the quantile information). The authors show that under some conditions on $\varepsilon$, the function $g(1-z)=\max\{\varepsilon,\sum_{j\in[m]}\{\max_{i \in [N]}\{r_{i,j}(1-z_i)\}\}$ is submodular. This implies, from Theorem \ref{thm:polymat} that the  polymatroid inequalities, for $t=\sum_{j\in[m]} t_j$, are sufficient to describe the convex hull of solutions to this set under the given conditions on $\varepsilon$. These are the so-called aggregated mixing inequalities for CCPs. 

}

\bibliographystyle{abbrvnat}
\bibliography{mybibfile,master,SIP-BIB}

\end{document}